\documentclass[aap,preprint]{imsart}

\RequirePackage[OT1]{fontenc}
\RequirePackage{amsthm,amsmath}
\RequirePackage[numbers]{natbib}
\RequirePackage[colorlinks,citecolor=blue,urlcolor=blue]{hyperref}

\startlocaldefs
\numberwithin{equation}{section}
\theoremstyle{plain}

\endlocaldefs

\usepackage{mathrsfs}
\usepackage{amssymb}
\usepackage{amsmath}
\usepackage{amsbsy}
\usepackage{epsfig}
\usepackage{enumerate}
\usepackage{bm}
\usepackage{color}
\usepackage{booktabs}

\newtheorem{example}{{\sc Example}}
\newtheorem{lemma}{Lemma}[section]
\newtheorem{theorem}{Theorem}[section]

\newtheorem{assumption}{Assumption}[section]
\newtheorem{corollary}{Corollary}[section]
\newtheorem{remark}{\sc Remark}[section]
\newtheorem{definition}{\sc Definition}[section]

\def\proclaim#1{\par \smallskip\noindent {\bf #1}\bgroup\it\ }
\def\endproclaim{\egroup\par\smallskip}

\def\pr{\textsf{P}} 
\def\ep{\textsf{E}} 

\def\underwiggle 1{
\ifmmode\setbox\TempBox=\hbox{$ 1$}\else\setbox\TempBox=\hbox{1}\fi
\setbox\TempBoxA=\hbox to \wd\TempBox{\hss\char'176\hss}
\rlap{\copy\TempBox}\smash{\lower9pt\hbox{\copy\TempBoxA}}}

\begin{document}

\begin{frontmatter}

\title{Convergence of randomized urn models with irreducible and reducible replacement policy\thanksref{T1}}
\runtitle{Convergence of randomized urn models }

\begin{aug}

\author{\fnms{Li-Xin} \snm{ZHANG}\ead[label=e1]{stazlx@zju.edu.cn}
\ead[label=u1,url]{http://www.math.zju.edu.cn/zlx}}

\address{Li-Xin ZHANG\\
School of Mathematical Sciences$~~~~~~~~~~~~~~~$ \\
Zhejiang University\\
Zheda Road, NO. 38~~~~~~ \\
Hangzhou, 310027~, P.R. China \\
\printead{e1}\\
\printead{u1}
}

\thankstext{T1}{Research supported by grants from the NSF of China
(No.11731012), the 973 Program
(No. 2015CB352302), Zhejiang Provincial Natural Science Foundation (No: LY17A010016)    and the Fundamental
Research Funds for the Central Universities. }

\runauthor{ L-X. Zhang }

\affiliation{Zhejiang University}
\end{aug}

\begin{abstract}

Generalized Friedman urn is one of the simplest and most useful models
considered in probability theory. Since Athreya and Ney (1972) showed the almost sure convergence of urn proportions in a randomized urn model with irreducible replacement matrix   under the $L\log L$ moment assumption, this  assumption has been regarded as the weakest moment assumption, but the necessary  has never been shown.  In this paper, we study the strong and weak convergence of generalized Friedman urns. It is  proved that, when the random replacement matrix is irreducible in probability,  the sufficient and necessary moment assumption for the almost sure  convergence of the urn proportions is   that the expectation of the  replacement matrix is finite, which is less stringent than the $L\log L$ moment assumption, and when the replacement is reducible, the $L\log L$ moment assumption is the weakest sufficient condition. The rate of convergence and  the strong and weak convergence of non-homogenous generalized Friedman urns are also derived.
\end{abstract}



\end{frontmatter}



\section{ Introduction}\label{sectionintroduction}
\setcounter{equation}{0}

Urn models have long been considered
powerful mathematical instruments in many areas, including the physical
sciences, biological sciences, social sciences and engineering [Johnson and
Kotz, 1977; Kotz and Balakrishnan, 1997].   The P\'olya urn
 model was originally proposed to
model the problem of contagious diseases [Eggenberger and P\'olya, 1923]. Since
then, there have been numerous generalizations and extensions. Among
them, the generalized Friedman urn (also named as
generalized P\'olya urn   in literature)   is the most popular one in the literature [see  Athreya and Karlin, 1968; Athreya and Ney, 1972; Higueras et al, 2003, 2006; Janson, 2004; etc.].
The generalized Friedman urn is also a papular model of response-adaptive randomization in clinical trial studies  [see Wei and Durham, 1978; Wei, 1979; Smythe, 1996; Bai and
Hu, 1999, 2005; Hu and Zhang, 2004a; Hu and Rosenberger, 2006; Zhang,
Hu and Cheung, 2006; Zhang et al, 2011; etc].      In the generalized Friedman urn,  the urn starts with the urn composition
$\bm Y_0=(Y_{0,1}, \cdots, Y_{0,d})\in \mathbb{R}_+^d\backslash \{\bm 0\}$. At the
stage $m$ ( $m=1,2,\cdots$), a ball is drawn from the urn with instant replacement. If the ball is of type $k$, then  an  additional random number $D_{k,q}(m)$  of balls of type
$q$, $q=1,\cdots, d$, are added to the urn.  After $n$ draws and
generations, the urn composition is denoted by the row vector $\bm
Y_n=(Y_{n,1}, \cdots, Y_{n,d})$, where $Y_{n,k}$ stands for the number
of balls of type $k$ in the urn after the $n$th draw.
This
relation can be written as the following recursive formula:
\begin{equation}\label{eq1.1}
 \bm Y_m = \bm Y_{m-1}+\bm X_m \bm D_m,
 \end{equation}
 where $\bm D_m=\big(D_{k,q}(m)\big)_{k,q=1}^d $, and $\bm X_m$ is the result
of the $m$th draw, distributed according to the urn composition at
the previous stage, i.e., if the $m$th draw is a type $k$ ball,
then the $k$th component of $\bm X_m$ is $1$ and other components
are $0$. The matrices $\bm D_m$'s are named as the   adding
rules or replacement matrices. We assume that $\bm D_m$, $m=1,2,\cdots, $ are independent and identically distributed (i.i.d.), and let
$\bm H=\ep[\bm D_m]=\big(\ep[D_{k,q}(m)]\big)_{k,q=1}^d $ when the expectations in the brackets are finite. $\bm H$ is said to be the generating matrix. In the original P\'olya urn model and many of its generalizations, $\bm D_n$ is a deterministic matrix. To distinguish whether $\bm D_n$ is deterministic or not, the urn model is also called a randomized urn model when $\bm D_n$  is random [see Bai and Hu, 2005].

When   $D_{k,q}(m+1)$, $Y_{m,k}$, $k,q=1,\ldots, d$, $m=0,2,\ldots$, all take non-negative integer values, Athreya  and Karlin  (1968) and Athreya and Ney (1972) studied the convergence of $\bm Y_n/n$ by embedding the urn process $\bm Y_n$ into a multi-type branching process. Let $\{\bm Z(t); t\ge 0\}$ be a $d$-type branching process for which (i) the life times of particles of all types are unit exponentials; and (ii) particles live and produce their offsprings  independently of each other, and of the past;   a $k$th type particle creates, on death, a random number $D_{k,k}+1$ of new particles of type $k$, and  a random number $D_{k,q}$ of new particles of type $q$, $q\ne k$, where the random vector $(D_{k,1},\ldots, D_{k,d})$ has the same distribution as $(D_{k,1}(m),\ldots, D_{k,d}(m))$ for each $k=1,\ldots, d$.  The matrix $\bm D=(D_{k,l})_{k,l=1}^d$ represents  the offspring producing rule.  Let $\{\tau_n; n=0,1,2,\ldots; \tau_0=0\}$ denote the split times for this process.  Athreya and Ney (1972) presented the link between the urn model and the branching process by the following embedding theorem, see Theorem 9.2 of Athreya and Ney (1972).
\begin{proclaim}{Theorem A} If  $\bm Z(0)=\bm Y_0$, then
the stochastic process $\{\bm Y_n;n=0,1,2,\ldots\}$ and  $\{\bm Z(\tau_n);n=0,1,2,\ldots\}$  are equivalent.
\end{proclaim}
 By this equivalence and the limit theorems for branching processes, Athreya and Ney (1972) showed the almost sure convergence of $\frac{\bm Y_n}{n}$ under the $L\log L$ moment assumption, namely,
\begin{equation}\label{LlogLcondition}
\ep[ D_{k,q} \log  D_{k,q} ]<\infty \;\; \text{ for all } k,q,
\end{equation}
and an assumption on $\bm H$ that $\bm H^l$ is a matrix with positive entries for some integer $l\ge 1$, see also Athreya and Ney (2004). This assumption on $\bm H$ implies that the nonnegative matrix $\bm H$ is irreducible (for the definition see the next section) and so the largest real part $\lambda_H$ of  all eigenvalues of $\bm H$ is a simple and positive eigenvalue and, associated with $\lambda_H$  the    nonnegative left eigenvector $\bm v=(v_1,\cdots,v_d)$ with $v_1+\cdots v_d=1$ is unique and positive. The limit of $\frac{\bm Y_n}{n}$ is just $\lambda_H \bm v$.
Because (\ref{LlogLcondition}) is also a necessary condition for $e^{-\lambda_H t}\bm Z(t)$ converging to a non-zero vector (c.f. Theorem 7.2 of Athreya and Ney, 1972),   for a long history it has been   expected that   (\ref{LlogLcondition}) is the weakest moment condition for studying the almost sure convergence of the generalized Friedman urn models. By the supermartingale method, Zhang (2012) proved a similar almost sure convergence of $\frac{\bm Y_n}{n}$  under (\ref{LlogLcondition}) when $D_{k,q}(m+1)$, $Y_{m,k}$, $k,q=1,\ldots, d$, $m=0,2,\ldots$,  take non-negative   real, not necessary integer, values.
Though the generalized Friedman urn model has been extended and studied in various ways  [c.f.   Bai
and Hu, 1999, 2005; Bena\"{\i}m, Schreiber  and Tarr\`{e}s, 2004; Higueras, et al, 2003, 2006; Janson, 2004; Laruelle  and Pag\'es, 2013; Zhang, Hu and Cheung 2006; etc], the almost sure convergence  is usually showed under moment conditions more stringent   than  (\ref{LlogLcondition}). For examples, by applying the theory of branching process, Janson (2004) obtained the almost sure convergence and central limit theorems under the second moment finite, namely, $\ep[\|\bm D_n\|^2]<\infty$. By applying the theory of matrices and the stochastic approximation (SA)
algorithm, respectively, Bai and Hu (2005) and Laruelle and Pag\'es (2013) showed the almost sure convergence  for the case of non-homogeneous replacement in which $\{\bm D_n\}$ may be not i.i.d.,  under the second moment finite, but the replacement is assumed to be balanced, namely, all the row sums of the (conditional) expectation of $\bm D_n$ are equal. For the general non-homogeneity case that the replacement may be unbalanced, Zhang (2016) proved the almost sure convergence by studying the stability of a stochastic approximation algorithm with a non-linear regression function, but the second moment finite is still assumed. The central limit theorems are also obtained in Bai and Hu (2005),  Laruelle and Pag\'es (2013) and Zhang (2016) under the $(2+\epsilon)$-th moment finite. It seems that the almost sure convergence $\frac{\bm Y_n}{n}$ has not been  proved under a moment condition less stringent  than (\ref{LlogLcondition}).    Whether (\ref{LlogLcondition}) is necessary  has also never been shown in literature.

 On the other hand, in the studies on the generalized  Friedman urn models with random replacement matrices $\bm D_n$s, the irreducibility of the mean matrix $\bm H$  is usually an essential condition.
  The randomized urn model with   replacement  reducible  is seldom studied in literature.
    The purpose of this paper is to find the sufficient and necessary moment assumption for the almost sure convergence of $\frac{\bm Y_n}{n}$ in a generalized Friedman urn model with irreducible replacement or reducible replacement. We will find that, when the replacement is irreducible, the sufficient and necessary moment assumption is not (\ref{LlogLcondition}), but that  the expectations $\ep[D_{k,q}(n)]$, $k,q=1,\ldots,d$, are finite; and when   the replacement is reducible, the condition (\ref{LlogLcondition}) is a sufficient moment condition that can not be weakened and the limit proportions are  random.   The results obtained in this paper give a full picture of the convergence for both the   irreducible and reducible  replacement cases.  The rate of convergence, and weak and strong convergence of urn models with  non-homogeneous replacement     are also studied. The main results are given in the next section. The results for urns with non-homogenous replacement are given in Section \ref{section3}. The proofs are stated in the last section.  We will apply the method of the stochastic approximation algorithm to show the convergence as in Zhang (2016).
     For the stochastic approximation algorithm, the Kushner-Clark  theorem [c.f. Kushner-Clark, 1978;  Kushner  and Yin, 2003; Duflo, 1997] is a usual tool to show the almost sure  convergence. But now the  equilibrium point of the regression of the stochastic approximation algorithm may be not unique and the related  ordinary
differential equation (ODE) may be not stable. We will use a  direct way instead of the  Kushner-Clark theorem to find the limit.  For considering the properties of the limit proportions in the reducible replacement case, we will apply the   supermartingale to show that when $k$ is in an irreducible class corresponding the largest real eigenvalue of $\ep[\bm D_n]$, the limit of $\frac{\bm Y_{n,k}}{n}$ is a positive random variable  and zero otherwise,  and derive a  conditional central limit theorem to show that the positive limit having no point probability mass.

\section{ Main results}\label{sectionmain}
\setcounter{equation}{0}
Before we state the results. We first need some more notations and assumptions. As in Zhang (2016),  to include various cases, we allow the numbers of balls to be non-integers and negative. For example, $D_{k,l}(n)<0$ means that $|D_{k,l}(n)|$ balls of type $l$ is removed from the urn when a ball of type of $k$ is drawn. We assume that a type of ball with a negative number  will never be selected and so the selection probabilities are
\begin{equation}\label{eqassign} \pr(X_{n,k}=1|\mathscr{F}_{n-1})=\frac{ Y_{n-1,k}^+}{\sum_{j=1}^d  Y_{n-1,j}^+}, \;\; k=1,\cdots, d.
\end{equation}
Here $Y_{n,k}^+=\max\{Y_{n,k},0\}$ is the positive part of $Y_{n,k}$, $\bm Y_n^+=(Y_{n,1}^+,\cdots, Y_{n,d}^+)$,   $\mathscr{F}_m$ is the history sigma-field generated by $\bm X_1,\ldots,\bm X_m$, $\bm Y_1,\ldots,\bm Y_m$, $\bm D_1$, $\ldots$, $\bm D_m$,    and $\frac{\bm 0}{0}$ is defined to be $(p_1,\cdots,p_d)$, which means that a $k$-type   ball  is selected with  probability $p_k$ when the urn has no balls with a positive number, $k=1,\ldots,d$. Here $p_1,\cdots,p_d$ are pre-specified   probabilities with $\sum_kp_k=1$.  In this general framework, the urn allows the negative and/or non-integer  number of balls, and the removal.
Write $\bm N_n=(N_{n,1},\cdots, N_{n,d})$, where
$N_{n,k}$ is the number of times that a type $k$ ball is drawn in the
first $n$ stages.
Obviously, $
 \bm N_n=\sum_{k=1}^n \bm X_n.
 $

In a general  branching process $\bm Z(t)$, the life times of particles may not have the same distribution. If the life times of particles of  type $k$ are    exponential  with parameter $\alpha_k>0$, $k=1,\cdots, d$, then the related urn processes is that with selection probabilities defined as
$$ \pr(X_{n,k}=1|\mathscr{F}_{n-1})=\frac{ \alpha_k Y_{n-1,k}}{\sum_{j=1}^d \alpha_j Y_{n-1,j}}, \;\; k=1,\cdots, d.
$$
Janson (2004) studied the properties of this kind of urn models. Under our framework, because the balls allow non-integer numbers, we can redefine the urn process  as $(\alpha_1Y_{n,1},\cdots, \alpha_dY_{n,k})$ with replacement matrices $\bm D_n diag(\alpha_1,\cdots,\alpha_d)$. The redefined urn process satisfies (\ref{eq1.1}) with (\ref{eqassign}) and generating matrix  $\bm Hdiag(\alpha_1,\cdots,\alpha_d)$.

For considering the asymptotic properties, we need assumptions on the replacement matrices.

\begin{assumption} \label{assumption4.1}  $\{\bm D_n\}$ is a sequence of independent and identically distributed random matrices.
\end{assumption}

 \begin{assumption} \label{assumption4.2} The expectations $H_{k,q}=\ep[D_{k,q}(n)]$, $k,q=1,\ldots, d$, are finite. Let $\bm H=(H_{q,k})_{k,q=1}^d$.  \end{assumption}

For a $d\times d$-matrix $\bm H$, we denote
$$\lambda_H=\max\{Re(\lambda): \lambda \text{ is an eigenvalue of } \bm H\} $$
to be the largest real part of its eigenvalues. If $H_{k,q}\ge 0$ for $q\ne k$, then there is a $\beta>0$ such that $(\beta I_d+\bm H)$ is a nonnegative matrix, where $\bm I_d$ is the $d\times d$ -identity matrix, and so  (i) $\lambda_H$ is an eigenvalue of $\bm H$,  (ii) if $\lambda\ne \lambda_H$ is an eigenvalue of $\bm H$, then $Re(\lambda)<\lambda_H$, and (iii) $\bm H$ has nonnegative  left eigenvectors and nonnegative right eigenvectors of $\bm H$ corresponding to $\lambda_H$. We let $\bm S_H=\{\bm v: \bm v\bm H=\lambda \bm v, \bm v\in  \bm \Delta^d\}$, where $\bm \Delta^d=\{\bm v: \sum_{k=1}^d v_k=1, v_k\ge 0, k=1,\ldots, d\}$, to denote the space of scaled nonnegative left eigenvectors corresponding to $\lambda_H$.

A  square matrix $\bm A$ is said to be reducible when there exists a permutation matrix $\bm P$ such that
 $$ \bm P^t\bm A\bm P=\begin{pmatrix}\bm A_{11} & \bm A_{12} \\ \bm 0 & \bm A_{22}\end{pmatrix},\; \text{ where } \bm A_{11} \text{ and } \bm A_{22} \text{ are both square}. $$
 Otherwise $\bm A$ is said to be an irreducible matrix. A $d\times d$ matrix $\bm A$ with nonnegative off-diagonal entries  is irreducible if and only if $\beta\bm I_d+\bm A\ge \bm 0$ and $(\beta\bm I_d+\bm H)^{d-1}>\bm 0$ for some $\beta>0$. Here and in the sequel, $\bm A\ge\bm 0$ and $\bm A>\bm 0$ mean that  all entries of $\bm A$ are nonnegative and positive, respectively.   By the Perron-Frobenius theory, if $\bm H$ with $H_{k,q}\ge 0$ ($k\ne q$)   is irreducible, then  (i)    $\lambda_H$  is a simple  eigenvalue of $\bm H$;
(ii) There exist an unique right eigenvector  $\bm v=(v_1,\cdots,v_d)$ and left
eigenvector $\bm u^{\rm t}=(u_1,\cdots, u_d)^{\rm t}$ corresponding to $\lambda_H$ such that   $\sum_k
v_k=\sum_k v_k u_k=1$ and $ v_k> 0, u_k>0$, $k=1,\cdots, K$, and so $\bm S_H=\{\bm v\}$.  In general, we assume the following assumption.

 \begin{assumption} \label{assumption4.3}

\begin{description}
  \item[\rm (a)] $H_{q,k}\ge 0$ for $q\ne k$;
  \item[\rm (b)] there exist a right
eigenvector $\bm u^{\rm t}=(u_1,\cdots, u_d)^{\rm t}$ corresponding to $\lambda_H$ such that  $u_k>0$, $k=1,\cdots, K$; and
  \item[\rm (c)]  $\lambda_H>0$.
\end{description}
\end{assumption}

\begin{theorem}\label{th2.1}  If Assumptions \ref{assumption4.1}-\ref{assumption4.3} are satisfied, then
\begin{equation} \label{eqth2.1.1}
\lim_{n\to \infty}dist\Big(\frac{\bm Y_n}{n},\lambda_H\bm S_H\Big)=\lim_{n\to \infty}dist\Big(\frac{\bm Y_n^+}{n},\lambda_H\bm S_H\Big)=0\;\; a.s.,
\end{equation}
\begin{equation} \label{eqth2.1.2}
\lim_{n\to \infty}dist\Big(\frac{\bm Y_n^+}{\sum_{j=1}^d Y_{n,j}^+}, \bm S_H\Big)=0\;\; a.s.
\end{equation}
and
\begin{equation} \label{eqth2.1.3}
\lim_{n\to \infty}dist\Big(\frac{\bm N_n}{n},\bm S_H\Big) \;\; a.s.
\end{equation}
Here the distance $dist(x,\bm S)$ between a point $\bm x$ and a set $\bm S$ is defined by $dist(\bm x,\bm S)=\inf\{\|\bm x-\bm y\|: \bm y\in \bm S\}$.
\end{theorem}

 The following theorem give the necessariness of the  expectations $\ep[D_{k,q}(n)]$s being finite when the limit proportions are positive.

\begin{theorem}\label{th2.2}
Suppose that Assumption \ref{assumption4.1} is satisfied  and $\ep[D_{q,k}(n)]>-\infty$ for all $q,k$.

(a) If there is a random vector $\bm V$ with  $ \pr(\bm V>\bm 0)>0$, such that
\begin{equation}\label{eqth2.2.1}
\frac{\bm Y_n}{n}\to  \bm V\;\; \text{in probability},
\end{equation}
  then Assumption \ref{assumption4.2} is satisfied.

(b) Suppose that there is a random vector $\bm V$ with $\sum_{k=1}^d V_k=1$ and $ \pr(\bm V_k>\bm 0)>0$, such that $\frac{Y_{n,k}}{\sum_{j=1}^d Y_{n,j}}\to V_k$ or $\frac{N_{n,k}}{n}\to V_k$   in probability for  $k=1,\cdots, d$. If the entries in one column of $\bm H=\ep[\bm D_n]$ are finite, then all entries of $\bm H$ are finite.
\end{theorem}

 Usually, as in Athreya  and Karlin  (1968) and Athreya and Ney (1972) etc,  the values $D_{k,q}(n)$s are assumed to be non-negative. Sometimes,  as in Janson (2004) and Laruelle  and Pag\'es (2013) etc, the drawn ball is allowed to be dropped, so the diagonal elements  $D_{k,k}(n)$s can take negative value $-1$ or $-c_k$ for some positive constant $c_k$. Such replacement matrices satisfy the following assumption.

\begin{assumption}\label{assumption4.4}  $D_{k,q}(n)\ge 0$ a.s. for all $k\ne q$, and $D_{k,k}(n)\ge -c_0 $ for some $c_0$ and all $k$.
\end{assumption}

A simple example of the replacement matrix satisfies Assumption \ref{assumption4.3} but not  Assumption \ref{assumption4.4} is that $\widetilde{D}_{k,q}(n)=D_{k,q}(n)+\epsilon_{k,q}(n)$, where $D_{k,q}(n)$s satisfy Assumption \ref{assumption4.4} and $\epsilon_{k,q}(n)$ are replacement errors  with mean zeros.

When Assumption \ref{assumption4.4} is satisfied, then for any $C>0$ the matrix $\ep[\bm D_n\wedge C]=:\big(\ep[ C\wedge D_{q,k}(n)]\big)_{q,k=1}^d$ satisfies Assumption \ref{assumption4.3}(a).
It is easily  that, when $\ep[\bm D_n\wedge C]$ is irreducible, then $\ep[\bm D_n\wedge l]$ is irreducible for all $l\ge C$, because $\ep[\bm D_n\wedge l]\ge \ep[\bm D_n\wedge C]$ if $l\ge C$.
The irreducibility   of the matrix $\ep[\bm D_n]$ or $\ep[\bm D_n\wedge C]$ only depends on the structure of what off-diagonal elements are nonzero. If define $\bm P_n=\big(p_{k,q}(n)\big)_{k,q=1}^d$  by $p_{k,q}(n)=\pr(D_{k,q}(n)\ne 0)$, then $\bm P_n$ will have the same structure of $\ep[\bm D_n\wedge C]$ when $C$ is large. So, when $\bm H$ is not finite, we may define the irreducibility of the replacement by the  irreducibility of $\ep[\bm D_n\wedge C]$  or $\bm P_n$.

\begin{definition} The replacement matrix $\bm D_n$ is said to be irreducible in probability when $\bm P_n$ is irreducible, where $\bm P_n=\big(p_{k,q}(n)\big)_{k,q=1}^d$  with $p_{k,q}(n)=\pr(D_{k,q}(n)\ne 0)$.

The replacement matrix $\bm D_n$ satisfying Assumption \ref{assumption4.4} is said to be irreducible in mean when  $\ep[\bm D_n\wedge C]$ is irreducible for some $C>0$.
\end{definition}
It is obvious that, under Assumption \ref{assumption4.4}, $\bm D_n$ is irreducible in probability if and only if it is irreducible in mean, and when $\ep[\bm D_n]$ is finite,  the irreducibility  in mean,  irreducibility  in probability and that $\ep[\bm D_n]$ is irreducible are equivalent. So, in general, when $\bm D_n$ is irreducible in probability, the replacement is said to be irreducible.

Now, we consider the case of irreducible replacement. The following corollary follows from Theorem \ref{th2.1} immediately.

\begin{corollary}\label{corollary2.1}  Suppose that  Assumptions \ref{assumption4.1}-\ref{assumption4.3} are satisfied. Further,  assume that $\bm H$  is irreducible. Then
\begin{equation} \label{eqcorollary2.1.1}
\lim_{n\to \infty}\frac{\bm Y_n}{n}=\lim_{n\to \infty}\frac{\bm Y_n^+}{n}= \lambda_H\bm v\;\; a.s.,
\end{equation}
\begin{equation} \label{eqcorollary2.1.2}
\lim_{n\to \infty}\frac{\bm Y_n^+}{\sum_{j=1}^d Y_{n,j}^+}=\bm v\;\; a.s.
\end{equation}
and
\begin{equation} \label{eqcorollary2.1.3}
\lim_{n\to \infty}\frac{\bm N_n}{n}= \bm v \;\; a.s.,
\end{equation}
where  $\bm v\in \bm\Delta^d$ is the unique solution of  the equation $\bm v\bm H=\lambda_H\bm v$.
\end{corollary}

 The following is the converse of Corollary \ref{corollary2.1}.
\begin{corollary}\label{corollary2.2}  Suppose that Assumptions \ref{assumption4.1}, \ref{assumption4.4} are satisfied, and
   $\bm D_n$ is irreducible in probability. If there is a random vector $\bm V$ with  $ \pr(\bm V\ge \bm 0, \bm V\ne \bm 0)>0$ such that (\ref{eqth2.2.1}) holds, then Assumptions  \ref{assumption4.2}, \ref{assumption4.3} are satisfied, and $\bm H$ is irreducible. Further, $\bm V=\lambda_H\bm v$ a.s., where $\bm v\in \bm\Delta^d$ is the unique solution of $\bm v\bm H=\lambda_H\bm v$.
 \end{corollary}

 From Corollaries \ref{corollary2.1} and \ref{corollary2.2}, we conclude that, in the case of irreducible replacement,  the mean replacement matrix being finite is the sufficient and necessary condition for the proportions $\frac{\bm Y_n}{n}$ to have a non-zero limit.

 In literature, the studies on urn models with reducible replacements are very few.  Gouet (1997) considered the case of fixed deterministic balanced, but not necessarily
irreducible, replacement matrix.   Abraham,   Dhersin  and   Ycart  (2007) considered a  special  urn scheme with reducible random replacements, in which, each time a ball is picked, another ball is added, and its type is chosen according to the transition probabilities of a reducible Markov chain. The vector of frequencies $\bm Y_n/n$ is shown  to converge almost surely to a random element of the set of stationary measures of the Markov chain.

 Our next theorem shows that for the  randomized urn model, when $\bm H$ is reducible, $\frac{\bm Y_n}{n}$ and $\frac{\bm N_n}{n}$  will also converge almost surely under $L\log L$ moment assumption. To state the result, we need more notations  for describing  the structure of the matrix $\bm H$.  Let $\nu_1$ be the index of the eigenvalue $\lambda_H$. Under  (a) and (b) in Assumption \ref{assumption4.3} ,
$\bm H$   has the following Jordan canonical form
\begin{equation}\label{eqJordanform}
  \bm H = \bm T diag(\lambda_H\bm I_{\nu_1}, \bm J_2,\cdots,\bm J_s)\bm T^{-1},
  \end{equation}
where
$$\bm J_t=\begin{pmatrix}
\lambda_t &1           &  0     &\ldots     & 0 \\
         0         & \lambda_t  &1       &\ldots     & 0 \\
         \vdots    & \ldots     & \ddots &\ddots     & \vdots \\
         0         & 0          &\ldots  & \lambda_t & 1 \\
         0         & 0          &0       & \ldots    &\lambda_t
\end{pmatrix}_{\nu_t\times\nu_t}, \;\; Re(\lambda_t)<\lambda_H. $$
Further, in (\ref{eqJordanform}), $\bm T$ and $\bm T^{-1}$ can be chosen such that $\bm v_1,\ldots, \bm v_{\nu_1}\in \bm S_H$ be the first $\nu_1$ rows of $\bm T^{-1}$, and $\bm u_j^t\ge 0$, $j=1,\ldots, \nu_1$, be the first $\nu_1$ columns of $\bm T$. Then $\bm H\bm u_j^t=\lambda_H\bm u_j^t$, $\bm v_j\bm u_j^t=1$, $\bm v_j\bm u_i^t=0$, $i\ne j$. $\bm S_H$ has the form $\{\bm v=\sum_{j=1}^{\nu_1} \beta_j \bm v_j: \sum_{j=1}^{\nu_1} \beta_j=1, \beta_j\ge 0, j=1,\cdots,\nu_1\}$. The set $\{1,\ldots,d\}$ can be divided to several irreducible classes $\mathcal{C}_j$, $j=1,\ldots, r$, such that (i) each principal submatrix $\bm H_{jj}=(H_{k,q};k,q\in \mathcal{C}_j\}$ of $\bm H$ on the class $\mathcal{C}_j$ is irreducible or  $(0)_{1\times 1}$; (ii)
$\lambda_{H_{jj}}=\lambda_H$, $j=1,\ldots,\nu_1$, and $\lambda_{H_{jj}}<\lambda_H$ for other $j$; (iii) the elements of   $\bm v_j$ in the class $\mathcal{C}_j$ are positive and zeros otherwise, $j=1,\cdots, \nu_1$; (iv) the elements of   $\bm u_j$ in the class $\mathcal{C}_j$ are positive and zeros in $\mathcal{C}_i$ for $i\ne j$, $j,i=1,\cdots,\nu_1$. Further, $\bm H$ has the following structure. If $q$ is in one of the classes $\mathcal{C}_i$, $i\le \nu_1$, and $k$ is not in the same class, then $H_{q,k}=0$.  For any $k$   in one of the classes of $\mathcal{C}_i$, $i> \nu_1$, there is class $\mathcal{C}_j$, $j=1,\cdots,\nu_1$, and a $q$ in it such that $H_{k,q}\ne 0$.

 Let $\bm U=:\sum_{j=1}^{\nu_1} \bm u_j^t\bm v_j$. Then $\bm U$  is a project onto the space
$\{\bm v:\bm v\bm H=\lambda_H\bm v\}=\{\bm v: \bm v=\sum_{j=1}^{\nu_1} \beta_j\bm v_j\}$.  It is obvious that $\bm H$ is irreducible if and only if $\nu_1=1$, and then $\bm S_H=\{\bm v\}$ has an unique point.

\begin{theorem}\label{th2.3}  Suppose that  Assumptions \ref{assumption4.1}-\ref{assumption4.3} are satisfied. Further,  assume
   $\ep\left[\|\bm D_n\|\log(\|\bm D_n\|)\right]<\infty$.
Then there exists a   random vector $\bm V$ which takes values in $\bm S_H$  such that
\begin{equation} \label{eqth2.3.1}
\frac{\bm Y_n}{n}\to \lambda_H\bm V\;\; a.s.,
\end{equation}
\begin{equation} \label{eqth2.3.2}
\frac{\bm Y_n^+}{\sum_{j=1}^d Y_{n,j}^+}\to \bm V\;\; a.s.
\end{equation}
and
\begin{equation} \label{eqth2.3.3}
\frac{\bm N_n}{n}\to \bm V \;\; a.s.
\end{equation}
If also $\bm D_n$s are non-negative  matrices and $Y_{0,k}>0$, $k=1,\cdots, d$, then the limit $\bm V\in \bm S_H$ has the form $\bm V=\sum_{j=1}^{\nu_1} \varpi_j\bm v_j$ with $\sum_{j=1}^{\nu_1 } \varpi_j=1$, $0<\varpi_j\le 1$, $j=1,\cdots, \nu_1$, with probability one, and further, if $\ep[\|\bm D_n\|^2]<\infty$, then each $\varpi_j$ has no point  probability mass  in $[0,1)$, namely, $\pr(\varpi_j=p)=0$ for any $p\in[0,1)$.
\end{theorem}

For the reducible case ($\nu_1>1$), from Theorem \ref{th2.3}, we find that if $\bm D_n$ is nonnegative and has second moment finite, then $\frac{Y_{n,k}}{\lambda_H n}$ and $\frac{N_{n,k}}{n}$ converge to a random variable $ V_k$ in $(0,1)$ when $k$ belongs to one of the  irreducible classes $\mathcal{C}_j$, $j=1,\cdots, \nu_1$, associated with the largest eigenvalue $\lambda_H$, and this random variable $V_k$ has no point mass.
When $k$ belongs to other classes, $\frac{Y_{n,k}}{\lambda_H n}$ and $\frac{N_{n,k}}{n}$ converge to zero. We conjuncture that the condition $\ep[\|\bm D\|^2]<\infty$ can be weakened to  $\ep\left[\|\bm D_n\|\log(\|\bm D_n\|)\right]<\infty$, and $\varpi_j$, $j=1,\ldots, \nu_1$ are continuous random variables having  densities with support $(0,1)$.

For the urn considered in Abraham,   Dhersin  and   Ycart  (2007), the replacement matrices $\bm D_n$ take zero or one entries with  $\pr(D_{k,q}(n)=1)=p_{k,q}$,   $\sum_{q=1}^d p_{k,q}=1$, and $\bm H=\ep[\bm D_n]=(p_{k,q})_{k,q=1}^d$  possibly reducible.
This model satisfies Assumptions \ref{assumption4.1}-\ref{assumption4.4} with $\bm H=\bm P$, $\lambda_H=1$ and $\bm u=(1,\cdots, 1)$.  Abraham, Dhersin  and   Ycart  (2007) characterized the limit probability distribution  as the solution to a fixed point problem. Examples showed that the limit probability distribution has density.

The next example shows that, the condition $\ep\left[\|\bm D_n\|\log(\|\bm D_n\|)\right]<\infty$ can not be weakened for the reducible case.
 \begin{example}
Suppose that $\{\bm D_n=diag(D_{1,1}(n),\cdots, D_{d,d}(n)); n\ge 1\}$  is a sequence of i.i.d. diagonal   and non-negative matrices. This is a reinforced
urn studied by  Zhang et al (2014). Suppose that there is a random vector $\bm V$ with  $ \pr(\bm V>\bm 0)>0$  such that (\ref{eqth2.2.1}) holds. Then by Theorem \ref{th2.2} (a), $m_k=\ep[D_{k,k}(1)]<\infty$ for all $k$.

It is obvious that $m_k\ne 0$, for otherwise we will have $D_{k,k}(n)=0$ a.s., and then $Y_{n,k}=Y_{0,k}$ a.s. for all $n$. By  Theorem 2.3 of Zhang et al (2014), $m_1=\cdots=m_d$. So,
Assumptions \ref{assumption4.1}-\ref{assumption4.3} are satisfied with $\lambda_H=m_1$ and $\bm u=(1,\cdots,1)$.

Further,  by  Theorem 2.3 of Zhang et al (2014),   if one of  $\ep\big[[D_{k,k}(1)]\log [D_{k,k}(1)]\big]$s is finite then all of them are finite.
\end{example}

From Theorems \ref{th2.1} and \ref{th2.3}, we get the following corollary on the branching process.

\begin{corollary} \label{corollary2.3} Suppose $\bm Z(t)$ is a branching process with nonnegative offspring producing rule $\bm D$ and life time parameters $\alpha_1,\cdots,\alpha_d$.  Let $\bm M=\big(\alpha_i\ep[D_{i,j}]\big)_{i,j=1}^d$ and    suppose the largest real part $\lambda_M$ of  the  eigenvalues of $\bm M$ is positive.  If $\bm M$ is irreducible, then as $t\to \infty$,
\begin{equation}\label{eqcor2.3.1} \frac{\bm Z(t)}{\sum_{j=1}^dZ_j(t)} \to \bm v \;\; a.s.,
\end{equation}
\begin{equation}\label{eqcor2.3.2}\frac{N_k(t)}{\sum_{j=1}^d N_j(t)}\to \frac{\alpha_kv_k}{\sum_{j=1}^d \alpha_j v_j}\;\; a.s. \;\; k=1,\cdots, d,
\end{equation}
where $N_j(t)$ is number of $k$th type particles died up  to time $t$, and $\bm v\in \bm\Delta^d$ is the unique solution of the equation $\bm v\bm M=\lambda_M\bm v$ and positive. If $\bm M$ is reducible  but has a positive right eigenvector corresponding to $\lambda_H$, and (\ref{LlogLcondition}) is satisfied, then there exists an random vector $\bm v$ taking values in the space $\bm S_M$ such that (\ref{eqcor2.3.1}) and (\ref{eqcor2.3.2}) holds.
\end{corollary}

The last  theorem in this section gives the rate of the convergence.
\begin{theorem}\label{th2.4} Suppose that  Assumptions \ref{assumption4.1}-\ref{assumption4.3} are satisfied. Further,  assume
   $\ep\left[\|\bm D_n\|^2 \right]<\infty$.
Then there exists a   random vector $\bm V$ which takes values in $\bm S_H$  such that
\begin{equation} \label{eqth2.4.1}
\left\|\frac{\bm Y_n}{n}-\lambda_H\bm V\right\|=O(b_n)\;\; a.s.,
\end{equation}
\begin{equation} \label{eqth2.4.2}
\left\|\frac{\bm Y_n^+}{\sum_{j=1}^d Y_{n,j}^+}- \bm V\right\|=O(b_n)\;\; a.s.
\end{equation}
and
\begin{equation} \label{eqth2.4.3}
\left\|\frac{\bm N_n}{n}\to \bm V\right\|=O(b_n) \;\; a.s.
\end{equation}
where $b_n$ is defined as follows. Let $\rho=\max\{Re(\lambda_j)/\lambda_H: j=2,\cdots,s\}$ be the ratio of the second largest real part of the eigenvalues of $\bm H$  to the largest one, and $\nu_{sec}=\max\{\nu_j: Re(\lambda_j)=\rho\lambda_H\}$ be the largest algebraic multiplicity of the eigenvalues with the second largest real part. Define $b_n$ by
\begin{equation}\label{eqth2.4.4}
b_n=\begin{cases}  n^{\rho-1}(\log  n)^{\nu_{sec}-1}, & \text{ if } \rho>1/2,\\
 n^{-1/2}(\log  n)^{\nu_{sec}-1}(\log\log n)^{1/2}, & \text{ if } \rho=1/2,\\
n^{-1/2}(\log\log n)^{1/2}, & \text{ if } \rho<1/2.
\end{cases}
\end{equation}
\end{theorem}

\bigskip
We will show Theorems \ref{th2.1},\ref{th2.3} and \ref{th2.4} in the last section after  establishing the results for the general  models with   non-homogenous replacements.
 Here we give the proofs of Theorem \ref{th2.2} and Corollary \ref{corollary2.2}.
To prove the results, we need a lemma at first.
\begin{lemma} \label{lem1} (a) Let $\xi_{n,k}=f(D_{k,1}(n),\cdots, D_{k,d}(n))$ be a Borel function of $(D_{k,1}(n),\cdots, D_{k,d}(n))$. Then on the event $\{N_{n,k}\to \infty\}$,
$$ \lim_{n\to \infty} \frac{\sum_{m=1}^n X_{m,k} \xi_{m,k}}{N_{n,k}} =\ep[\xi_{1,k}] \;\; a.s.  \;\; \text{ if }  \ep[|\xi_{1,k}|]<\infty, $$
and
$$  \frac{\big|\sum_{m=1}^n X_{m,k} (\xi_{m,k}-\ep[\xi_{1,k}])\big|}{\sqrt{ N_{n,k}\log\log N_{n,k}}} =O (1) \;  a.s.  \;  \text{ if }  \ep[|\xi_{1,k}|^2]<\infty. $$
(b) Let $\Delta \bm M_{n,1}= \bm X_n-\ep[\bm X_n|\mathscr{F}_{n-1}]$.  Then
$$ \|\bm M_{n,1}\|  =O\Big(\sqrt{n\log\log n}\Big)\;  a.s., \; \; \max_{m\le n}\|\bm M_{n,1}\|=O(\sqrt{n})\;
\text{in probability};$$

(c) Let   $\Delta \bm M_{n,2}=\bm X_n(\bm D_n-\ep[ \bm D_n])$. Under Assumptions \ref{assumption4.1} and \ref{assumption4.2},
$$ \frac{\bm M_{n,2}}{n}\to \bm 0\;\; a.s. $$
Further, if $\ep[\|\bm D_m\|^2]<\infty$, then
$$\|\bm M_{n,2}\|=O\Big(\sqrt{n\log\log n}\Big)\;  a.s.  $$
\end{lemma}

{\bf Proof.} The proof of (a) can be found in Hu and Zhang (2004b) (c.f. their Lemma A.4).  (b) is obvious because $\Delta \bm M_{n,1}$ is a sequence of bounded martingale differences.   (c) is a direct conclusion of (a). $\Box$.

\bigskip

{\bf Proof of Theorem \ref{th2.2}.}   For a vector $\bm \theta=(\theta_1,\cdots,\theta_d)$, we write $\alpha \big(\bm \theta\big)=\sum_{j=1}^d \theta_j$.

(a) Suppose (\ref{eqth2.2.1}) holds. Then on the event $\{\bm V>\bm 0\}$, $\frac{\bm Y_n}{\alpha(\bm Y_n)}\to \frac{\bm V}{\alpha(\bm V)}$ in probability, and so $\frac{\bm Y_n^+}{\alpha(\bm Y_n^+)}\to \frac{\bm V}{\alpha(\bm V)}$ in probability,
\begin{align}\label{eqproofth2.2.1}
  \frac{N_{n,k}}{n}=&\frac{1}{n}\sum_{m=1}^n \frac{Y_{m-1,k}^+}{\alpha(\bm Y_{m-1}^+)}+\frac{1}{n}\sum_{m=1}^n (X_{m,k}-\ep[X_{m,k}|\mathscr{F}_{m-1}])\nonumber\\
  =&\frac{1}{n}\sum_{m=1}^n \frac{Y_{m-1,k}^+}{\alpha(\bm Y_{m-1}^+)}+o(1) \;\; a.s.  \\
 \to &\frac{V_k}{\alpha(\bm V)}\;\; \text{in probability}, \nonumber
 \end{align}
$k=1,\cdots, d$,  by Lemma \ref{lem1}(b).  Note for each constant $M>0$,
\begin{align*}
&\frac{Y_{n,k}}{n}  =\frac{\sum_{m=1}^n \sum_{q=1}^d X_{m,q}D_{q,k}(m)}{n}\\
&\ge  \frac{\sum_{m=1}^n X_{m,l}D_{l,q}^+(m)}{n}-\sum_{q=1}^d\frac{\sum_{m=1}^n X_{m,k}D_{q,k}^-(m)}{n}\\
&\ge   \frac{  \sum_{m=1}^n X_{m,q} \big(M\wedge D_{q,k}^+(m)\big)}{n}-\sum_{q=1}^d\frac{\sum_{m=1}^n X_{m,k}\ep[D_{q,k}^-]}{n}+o(1)\;\; a.s.\\
&\ge  \frac{N_{n,q}}{n} \ep[M\wedge D_{q,k}^+]-\min_{q,k}\ep[D_{q,k}^-]+o(1) \;\; a.s. \\
& \to  \frac{ V_q}{\alpha(\bm V)} \ep[M\wedge D_{q,k}^+]-\min_{q,k}\ep[D_{q,k}^-]\;\; \text{ in probability }
\end{align*}
 on the event $\{\bm V>\bm 0\}$,  by Lemma \ref{lem1} (a). If $\ep[ D_{q,k}]=\infty$, then by letting $M\to \infty$ we conclude that $\frac{Y_{n,k}}{n}\to +\infty$ in probability on the event  $\{\bm V>\bm 0\}$ which contradicts the assumption (\ref{eqth2.2.1}). So, Assumption \ref{assumption4.2} is satisfied.

\medskip
(b)  Note that, if $\frac{Y_{n,k}}{\alpha(\bm Y_n)}\to V_k$ in probability, then on the event $\{\bm V>\bm 0\}$, $\frac{Y_{n,k}^+}{\alpha(\bm Y_n^+)}\to V_k$ in probability and
$ \frac{N_{n,k}}{n}  \to V_k $
in probability  by (\ref{eqproofth2.2.1}). So, it is sufficient to consider the case that $\frac{N_{n,j}}{n}\to V_j$ in probability for all $j=1,\cdots, d$.

Suppose that  the entries in $k$th column of $\bm H$ are finite, and there is an element of $\bm H$, say $H_{l,q}$, such that $H_{l,q}=\ep[D_{l,q}(m)]=\infty$.  Then,
\begin{align*}
\frac{|Y_{n,k}|}{n}\le &\frac{\sum_{m=1}^n \sum_{j=1}^d X_{m,j}|D_{j,k}(m)|}{n}\\
=&\sum_{j=1}^d \frac{N_{n,j}}{n} \ep[|D_{j,k}|]+o(1)  \;\; a.s.\\
  \to & \sum_{j=1}^d V_j \ep[|D_{j,k}|]\;\;\text{ in probability},
\end{align*}
by Lemma \ref{lem1} (a) and the fact that $\frac{N_{n,j}}{n}\to V_j$  in probability, and
$$ \frac{Y_{n,q}}{n}  \to \infty\text{ in probability on the event } \{\bm V>\bm 0\}. $$
It follows that
$$ \frac{Y_{n,k}^+}{\alpha(\bm Y_n^+)} \le \frac{|Y_{n,k}|}{Y_{n,q}}  \to 0 \;\; \text{ in probability on the event } \{\bm V>\bm 0\}. $$
Hence, by (\ref{eqproofth2.2.1}) $\frac{N_{n,k}}{n}\to 0$ in probability on the event $\{\bm V>\bm 0\}$, which  contradicts the assumption that $\frac{N_{n,k}}{n}\to V_k$ in probability. So, we conclude that all entries of $\bm H$ are finite.  $\Box$

\medskip
{\bf Proof of Corollary \ref{corollary2.2}.}
  Under Assumption \ref{assumption4.4}, $Y_{n,k}$ will not decrease when a type $q$ ball is drawn $q\ne k$ and it will not decrease when it becomes negative. So, $Y_{n,k}\ge    -c_0$. It follows that $\frac{Y_{n,k}^-}{n}\to 0$ a.s., $k=1,\cdots, d$.   Suppose (\ref{eqth2.2.1}) holds. Then, one the event $\{\bm V\ge\bm 0, \bm V\ne \bm 0\}$, $\frac{\bm Y_n}{\alpha(\bm Y_n)}\to \frac{\bm V}{\alpha(\bm V)}$ in probability,   $\frac{\bm Y_n^+}{\alpha(\bm Y_n^+)}\to \frac{\bm V}{\alpha(\bm V)}$ in probability, and so
$$ \frac{\bm N_n}{n}\to \frac{\bm V}{\alpha(\bm V)}\;\; \text{ in probability}, $$
by (\ref{eqproofth2.2.1}).  It follows that, one the event $\{\bm V\ge\bm 0, \bm V\ne \bm 0\}$,
\begin{align*}
\frac{\bm Y_n}{n}=& \frac{1}{n}\sum_{m=1}^n \bm X_m\bm D_m\ge \frac{1}{n}\sum_{m=1}^n \bm X_m\big(\bm D_m\wedge M)\\
=&\frac{\bm N_n}{n}\ep[\bm D_1\wedge M]+o_P(1)\to \frac{\bm V}{\alpha(\bm V)}\ep[\bm D_1\wedge M]
\text{ in probability},
\end{align*}
by Lemma \ref{lem1} (a). It follows that
$  \alpha(\bm V)\bm V\ge \bm V \ep[\bm D_1\wedge M]$  on the event $\{\bm V\ge\bm 0, \bm V\ne \bm 0\}$.
By the assumption that $\bm D_1$ is irreducible in probability which is equivalent  to that it is irreducible in mean, there exist $M>0$ and $\beta>0$ such that $\beta\bm I_d+ \ep[\bm D_1\wedge M]\ge \bm 0$ and $(\beta\bm I_d+ \ep[\bm D_1\wedge M])^{d-1}>\bm 0$. Hence
$(\beta+\alpha(\bm V))^d\bm V  \ge \bm V(\beta\bm I_d+ \ep[\bm D_1\wedge M])^{d-1}>\bm 0$ on the event $\{\bm V\ge\bm 0, \bm V\ne \bm 0\}$.  So, $\pr(\bm V>\bm 0)=\pr(\bm V\ge 0,\bm V\ne 0)>0$.
Therefore, by Theorem \ref{th2.2} (a) which has been proved, $\bm H=\ep[\bm D_n]$ is finite. Then $\bm H$ is also irreducible and so, its   unique left eigenvector $\bm v\in \bm \Delta^d$  corresponding to $\lambda_H$ is positive. By  Lemma \ref{lem1} (c) again,
$$ \frac{\bm Y_n}{n}=\frac{\bm N_n}{n}\bm H +\frac{\bm M_{n,2}}{n} =\frac{\bm N_n}{n}\bm H +o(1) \;\; \text{ in probability}, $$
which implies that $   \bm V= \frac{\bm V}{\alpha(\bm V)} \bm H$  on the event $\{\bm V\ge\bm 0, \bm V\ne \bm 0\}$. So,  $\frac{\bm V}{\alpha(\bm V)}=\bm v$ and $\lambda_H=\alpha(\bm V)>0$ on the event $\{\bm V\ge\bm 0, \bm V\ne \bm 0\}$. Hence, Assumptions \ref{assumption4.1}-\ref{assumption4.3} are satisfied. Finally, combing (\ref{eqcorollary2.1.1}) and (\ref{eqth2.2.1}) yields $\bm V=\lambda_H\bm v$ a.s.
  $\Box$

\section{ Non-homogenous replacement}\label{section3}
\setcounter{equation}{0}

In this section, we consider the case that $\{\bm D_n\}$ are not i.i.d. matrices. Let $\bm H_n=\ep[\bm D_n|\mathscr{F}_{n-1}]=\big(H_{q,k}(n)\big)_{q,k=1}^d$.
\begin{theorem}\label{th3.1}  Suppose that there exists a random matrix $\bm H$ such that
\begin{equation}\label{th3.1.1} \sum_{m=1}^n \|\bm H_n-\bm H\|=o(n) \;\; a.s.,
\end{equation}
and with probability one, $H(\omega)$ satisfies Assumption \ref{assumption4.3} in which $\lambda_H$, $\bm S_H$, $\bm u$ may depend on $\omega$.
 Assume
\begin{equation}\label{eqSLLNM} \|\bm M_{n,2}\|=o(n)\;\; a.s.
\end{equation}
 Then  (\ref{eqth2.1.1})-(\ref{eqth2.1.3}) hold. Further, if  $\bm H$ is irreducible almost surely, then $\bm S_H$ has   an unique point $\bm v$  and so (\ref{eqcorollary2.1.1})-(\ref{eqcorollary2.1.3}) hold.

The condition (\ref{eqSLLNM})  is satisfied if, for all $q,k$,  there is a sequence $c_m=c_{m,q,k}>0$ such that
\begin{equation}\label{th3.1.2}
\sum_{m=1}^{\infty}\pr\left(|D_{q,k}(m)|\ge  c_m \big|\mathscr{F}_{m-1}\right)<\infty \;\; a.s.,
\end{equation}
\begin{equation}\label{th3.1.3}
\sum_{m=1}^n \ep\left[|D_{q,k}(m)|I\{|D_{q,k}(m)|\ge    c_m\}\big|\mathscr{F}_{m-1}\right]=o(n)\;\; a.s.
\end{equation}
and
\begin{equation}\label{th3.1.4}
\sum_{m=1}^{\infty}\frac{ \ep\big[ \widetilde{D}_{q,k}^2(m)  \big|\mathscr{F}_{m-1}\big]}{m^2}<\infty\;\; a.s.,
\end{equation}
where $\widetilde{D}_{q,k}(m)=D_{q,k}(m) I\{|D_{q,k}(m)|\le    c_m\}-\ep\left[D_{q,k}(m) I\{|D_{q,k}(m)|\le    c_m\}|\mathscr{F}_{m-1}\right]$.
\end{theorem}

\begin{remark}\label{remark3.1} It is easily seen that (\ref{th3.1.2})-(\ref{th3.1.4}) with ($c_m=m$) are implied by
\begin{equation}\label{eqrk3.1.1}
  \begin{matrix}\text{ either }   \sup_m\ep\left[\|\bm D_m\|\log^{1+\epsilon}(\|\bm D_m\|)\big|\mathscr{F}_{m-1}\right]<\infty\; a.s. \\
   \text{ or }
\sup_m\ep[\|\bm D_m\|\log^{1+\epsilon}(\|\bm D_m\|)]<\infty
 \end{matrix}
  \end{equation}
  for some $\epsilon>0$.
  When $\{\bm D_m\}$ is a sequence of i.i.d. matrices, then   $\ep[|D_{k,q}(m)|]<\infty$ implies (\ref{th3.1.1}) and (\ref{th3.1.2})-(\ref{th3.1.4}) with $c_m=m$.
 \end{remark}
 For the non-homogeneity  case with $\bm H$ irreducible, Bai and Hu (2005) proved (\ref{eqcorollary2.1.1})-(\ref{eqcorollary2.1.3}) under the conditions that
\begin{equation}\label{condBaiHu1}
\bm H_m\ge \bm 0  \text{ and } \bm H_m \bm 1^{\prime}=c \bm 1^{\prime} \;\; a.s.,
\end{equation}
\begin{equation}\label{th3.3.1} \sum_{m=1}^{\infty}\frac{ \|\bm H_n-\bm H\|}{m}<\infty \;\; a.s.,
\end{equation}
\begin{equation}\label{condBaiHu2}
\sup_m \ep\left[\|\bm D_m\|^{2+\epsilon}|\mathscr{F}_{m-1}\right]<\infty \;\;a.s. \text{ for some } \epsilon>0.
\end{equation}
 By applying a  stochastic approximation algorithm, Laruelle and Pag\'es (2013) relaxed the conditions (\ref{th3.3.1}) and (\ref{condBaiHu2})  to $\bm H_n\to \bm H$ a.s. and the second moment finite that $\sup_m \ep[\|\bm D_m\|^2|\mathscr{F}_{m-1}]<\infty$ which is still more stringent than (\ref{eqrk3.1.1}).   The condition (\ref{condBaiHu1}) means that the updating of the urn at each stage is balanced.  Zhang (2016) removed this condition and proved (\ref{eqcorollary2.1.1})-(\ref{eqcorollary2.1.3}) under (\ref{th3.1.1}) and a moment condition that $\sum_{m=1}^n\ep[D_{q,k}^2(m)]=O(n)$.  This  moment condition  implies (\ref{th3.1.2})-(\ref{th3.1.4}) with $c_m=\infty$.   Zhang (2012) proved (\ref{eqcorollary2.1.1})-(\ref{eqcorollary2.1.3}) under (\ref{th3.1.1}) and (\ref{eqrk3.1.1}),  but assumed that the
entries of the replacement matrices $\bm D_n$ are non-negative.

It is easily shown that the condition (\ref{th3.3.1}) implies (\ref{th3.1.1}). Under this condition and moment assumptions a little  stricter than those in Theorem \ref{th3.1}, the urn proportions also converge almost surely to a random point in $\lambda_H\bm S_H$ even $\bm H$ is reducible.

\begin{theorem}\label{th3.3}  Suppose that there exists a random matrix $\bm H$ satisfying (\ref{th3.3.1}),
and with probability one, $H(\omega)$ satisfies Assumption \ref{assumption4.3} in which $\lambda_H$, $\bm S_H$, $\bm u$ may depend on $\omega$. Assume
\begin{equation}\label{eqConvergenceM}  \sum_{m=1}^{\infty} \frac{\max\limits_{i\le m}\|\bm M_{i,2}\|}{m^2}<\infty \;\;a.s.
\end{equation}
  Then     there exists a   random vector $\bm V$ taking values in $\bm S_H$  such that
(\ref{eqth2.3.1})-(\ref{eqth2.3.3}) hold.

The condition (\ref{eqConvergenceM}) is satisfied if,  for all $q,k$,  there is a sequence $c_m=c_{m,q,k}>0$ such that
\begin{equation}\label{th3.3.3}
\sum_{m=1}^{\infty}\frac{\ep\left[|D_{q,k}(m)|I\{|D_{q,k}(m)|\ge    c_m\}\big|\mathscr{F}_{m-1}\right]}{m}<\infty\;\; a.s.
\end{equation}
and
\begin{equation}\label{th3.3.4}
\sum_{m=1}^{\infty}\frac{(\log m)^{1+\epsilon} \ep\big[ \widetilde{D}_{q,k}^2(m)  \big|\mathscr{F}_{m-1}\big]}{m^2}<\infty\;a.s.\; \text{ for some } \epsilon>0,
\end{equation}
where $\widetilde{D}_{q,k}(m)=D_{q,k}(m) I\{|D_{q,k}(m)|\le    c_m\}-\ep\left[D_{q,k}(m) I\{|D_{q,k}(m)|\le    c_m\}|\mathscr{F}_{m-1}\right]$.
\end{theorem}

\begin{remark}
\label{remark3.2} It is easily seen that (\ref{eqrk3.1.1})  implies (\ref{th3.3.3}) and (\ref{th3.3.4})  with $c_m=m/(\log m)^{1+\epsilon}$.
When $\{\bm D_n;n\ge 1\}$ are i.i.d.,  the condition that $\ep\big[\|\bm D_n\|\log(\|\bm D_n\|)\big]<\infty$ implies (\ref{th3.3.3}) and (\ref{th3.3.4})   with $c_m=m/(\log m)^{1+\epsilon}$.
\end{remark}

In fact,  write $c_m=m/(\log m)^{1+\epsilon}$ and $\sigma_i^2=\ep[\widetilde{D}_{q,k}(i)^2|\mathscr{F}_{i-1}]$. Then
\begin{align*}
&\text{the left hand of (\ref{th3.3.3})} \\
\le &\sum_{m=1}^{\infty}\frac{1}{m} \frac{1}{(\log c_m)^{1+\epsilon}}\ep\big[|D_{q,k}(m)|\log^{1+\epsilon} (|D_{q,k}(m)|)|\mathscr{F}_{m-1}\big]\\
& \le     \sum_{m=1}^{\infty}\frac{1}{m} \frac{1}{(\log m)^{1+\epsilon}}\ep\big[|D_{q,k}(m)|\log^{1+\epsilon} (|D_{q,k}(m)|)|\mathscr{F}_{m-1}\big] \\
 & \le    c \sup_m\ep\big[|D_{q,k}(m)|\log^{1+\epsilon} (|D_{q,k}(m)|)|\mathscr{F}_{m-1}\big],
\end{align*}
and
\begin{align*}
&\text{the left hand of (\ref{th3.3.4})}=  \sum_{m=1}^{\infty}\frac{\sigma_m^2(\log m)^{1+\epsilon/2}}{m^2}\\
\le  &   \sum_{m=1}^{\infty}\frac{(\log m)^{1+\epsilon/2} \ep\big[|D_{q,k}(m)|^2I\{|D_{q,k}(m)|\le m/(\log m)^{1+\epsilon}\}\big|\mathscr{F}_{m-1}\big]}{m^2}\\
\le &c\sum_{m=1}^{\infty} \frac{m/(\log m)^{1+\epsilon/2}\ep\big[|D_{q,k}(m)|\log (|D_{q,k}(m)|)|\mathscr{F}_{m-1}\big]}{m^2}\\
\le &c\sup_m\ep\big[|D_{q,k}(m)|\log (|D_{q,k}(m)|)|\mathscr{F}_{m-1}\big],
\end{align*}
which both are finite when the first condition in (\ref{eqrk3.1.1}) is satisfied. When the second condition in (\ref{eqrk3.1.1}) is satisfied,   taking the expectation instead of the conditional expectation get the same conclusion.
If $\{\bm D_m\}$ are i.i.d., then
\begin{align*}
&\text{left hand of (\ref{th3.3.3})} = \sum_{m=1}^{\infty}\frac{1}{m} \ep\big[|D_{q,k}(1)|I\{|D_{q,k}(1)|\ge c_m\}\big]\\
 = &\ep\left[ |D_{q,k}(1)| \sum_{m=1}^{\infty}\frac{1}{m}  I\{|D_{q,k}(1)|\ge c_m\}\right]
 \le   c \ep[|D_{q,k}(1)|\log(|D_{q,k}(1)|)].
\end{align*}

The next theorem gives the rate of convergence.

\begin{theorem}\label{th3.4}  Suppose that there exists a random matrix $\bm H$ satisfying
\begin{equation}\label{eqth3.4.1}
\sum_{m=1}^n\|\bm H_m-\bm H\|=O(\sqrt{n\log\log n}) \;\; a.s.,
\end{equation}
and with probability one, $H(\omega)$ satisfies Assumption \ref{assumption4.3} in which $\lambda_H$, $\bm S_H$, $\bm u$ may depend on $\omega$. Assume
\begin{equation}\label{eqth3.4.2}   \|\bm M_{n,2}\|=O(\sqrt{n\log\log n})\;\;a.s.
\end{equation}
  Then     there exists a   random vector $\bm V$ taking values in $\bm S_H$  such that
(\ref{eqth2.4.1})-(\ref{eqth2.4.3}) hold, where $ b_n$ by defined as in (\ref{eqth2.4.4}), but now $\rho$ and $\nu_{\sec}$ are random variables.

The condition (\ref{eqth3.4.2}) is satisfied if (\ref{condBaiHu2}) is satisfied.
\end{theorem}

 The condition  (\ref{th3.1.1}) means that the conditional expectation $\bm H_n$ of the random replacement matrix $\bm D_n$ will converge almost surely to $\bm H$ in  Cesaro sense. If it    does not holds, then we may have no result on the convergence even high moments are assumed   finite. Recently,  Gangopadhyay and   Maulik (2017) showed that when the almost sure convergence in (\ref{th3.1.1}) is replaced by the convergence in probability, then $\frac{\bm Y_n}{n}$  converges in probability    under the conditions that $\bm H$ is irreducible, $\{\bm D_n; n\ge 1\}$ is   uniformly integrable, and $\bm Y_0$ has first moment finite.     The last theorem of this paper  gives a general result on the weak convergence of $\frac{\bm Y_n}{n}$.

\begin{theorem}\label{th3.5}  Suppose that there exists a random matrix $\bm H$ such that
\begin{equation}\label{th3.5.1} \sum_{m=1}^n \|\bm H_n-\bm H\|=o(n) \;\; \text{ in probability},
\end{equation}
and with probability one, $H(\omega)$ satisfies Assumption \ref{assumption4.3} in which $\lambda_H$, $\bm S_H$, $\bm u$ may depend on $\omega$. Assume
\begin{equation}\label{eqWLLNM}
\max_{m\le n}\|\bm M_{m,2}\|=o(n)\; \text{ in probability}.
\end{equation}
 Then, for any $T>0$,
\begin{align}\label{th3.5.4}
\nonumber
& \lim_{n\to \infty}\max_{n\le m\le nT}dist\Big(\frac{\bm Y_m}{m},\lambda_H\bm S_H\Big)\\
=&\lim_{n\to \infty}\max_{n\le m\le nT} dist\Big(\frac{\bm Y_m^+}{m},\lambda_H\bm S_H\Big)=0,
\end{align}
\begin{equation}\label{th3.5.5}
\lim_{n\to \infty}\max_{n\le m\le nT}dist\Big(\frac{\bm Y_m^+}{\alpha(\bm Y_m^+)}, \bm S_H\Big)=0
\end{equation}
and
\begin{equation}\label{th3.5.6}
\lim_{n\to \infty}\max_{n\le m\le nT}dist\Big(\frac{\bm N_m}{m}, \bm S_H\Big)=0
\end{equation}
in probability. In particular, if $\bm H$ is irreducible almost surely, then $\bm S_H$ has   an unique random point $\bm v$  and so
\begin{equation}\label{th3.5.7}
\lim_{n\to \infty}\frac{\bm Y_n}{n}= \lim_{n\to \infty}\frac{\bm Y_n^+}{n}=\lambda_H\bm v\;\; \text{ in probability},
\end{equation}
\begin{equation}\label{th3.5.8}
\lim_{n\to \infty}\frac{\bm Y_n^+}{\sum_{j=1}^d Y_{n,j}^+}= \bm v\;\; \text{ in probability}
\end{equation}
and
\begin{equation}\label{th3.5.9}
\lim_{n\to \infty}\frac{\bm N_n}{n}= \bm v\;\; \text{ in probability}.
\end{equation}

The condition (\ref{eqWLLNM}) is satisfied, if
\begin{equation}\label{th3.5.2}
\sum_{m=1}^n\ep\left[\|\bm D_m\|  \big|\mathscr{F}_{m-1}\right]=O(n) \;\; \text{ in probability}
\end{equation}
and for any $\epsilon>0$
\begin{equation}\label{th3.5.3}
\sum_{m=1}^n \ep\left[\|\bm D_m\|I\{\|\bm D_m\|\ge  \epsilon  n\}\big|\mathscr{F}_{m-1}\right]=o(n)\;\;\text{ in probability}.
\end{equation}
\end{theorem}

\begin{remark} If $D_{q,k}(n)\ge -c_0$ for some $c_0$, then (\ref{th3.5.2}) is implied by (\ref{th3.5.1}) and can be removed. The uniform integrability of $\{\bm D_n\}$ implies   (\ref{th3.5.3}).
\end{remark}
\begin{remark}
In Theorem \ref{th3.5}, (\ref{th3.5.4})-(\ref{th3.5.6}) give the results for the case of the reducible replacement, but we don't know whether there is a $\bm v$ such that (\ref{th3.5.7})-(\ref{th3.5.9}) hold in probability.
\end{remark}

\section{ Proofs}\label{sectionProof}
\setcounter{equation}{0}
 We apply the stochastic approximation algorithm method  as we did in Zhang (2016).   Note that $\bm Y_{n+1}=\bm Y_n+\bm X_{n+1}\bm D_{n+1}$,  $\ep[\bm X_{n+1}|\mathscr{F}_n]=\frac{\bm Y_n^+}{\alpha(\bm Y_n^+)}$ and $\ep[\bm D_{n+1}|\mathscr{F}_n]=\bm H_{n+1}$. Recall $\Delta \bm M_{n,1}= \bm X_n-\ep[\bm X_n|\mathscr{F}_{n-1}]$ and  $\Delta \bm M_{n,2}=\bm X_n(\bm D_n-\ep[ \bm D_n|\mathscr{F}_{n-1}])$.  We have
\begin{align*}
 \bm Y_{n+1}
 =&\bm Y_n+ \frac{\bm Y_n^+}{\alpha(\bm Y_n^+)} \bm H+\Delta \bm M_{n+1,1}\bm H +\Delta \bm M_{n+1,2}
 +\bm X_{n+1}\left( \bm H_{n+1} -\bm H\right).
  \end{align*}
  Write $\bm \theta_n=\frac{\bm Y_n^+}{n}$ and
\begin{align*}
\bm r_n = & (\Delta\bm M_{n,1}\bm H +  \Delta\bm M_{n,2}) +\bm X_{n+1}\left(\bm H_n-\bm H\right)+(\Delta\bm Y_n^+-\Delta\bm Y_n),\\
\bm s_n = & \sum_{m=1}^n \bm r_m= \bm M_{n,1}\bm H +  \bm M_{n,2} +\sum_{m=1}^n\bm X_m\left(\bm H_{m-1}-\bm H\right) +(\bm Y_n^+-\bm Y_n)
\end{align*}
and
\begin{equation}\label{eqregressF} \bm h(\bm \theta)=\bm\theta\Big(\bm I_d-\frac{\bm H}{\alpha(\bm \theta)}\Big).
\end{equation}
where $\alpha(\bm \theta)=\sum_{k=1}^d \theta_k$.  Then $\bm \theta_n$  satisfies the stochastic approximation  algorithm:
\begin{equation}\label{eqSA} \bm \theta_{n+1}=\bm\theta_n-\frac{\bm h(\bm\theta_n)}{n+1}+\frac{\bm r_{n+1}}{n+1}.
\end{equation}
Its related ordinary
differential equation (ODE) is
\begin{equation}\label{eqproofODE} \frac{d}{dt}\bm\theta(t)=-\bm h\big(\bm\theta(t)\big).
\end{equation}
It is obvious that every point $\lambda_H\bm v$  in $\lambda_H \bm S_H$ is an equilibrium point   of $\bm h(\bm\theta)$, namely, $\bm h(\lambda_H\bm v)=\bm 0$. When $\bm H$ is irreducible, Zhang (2016) show that $\bm\theta_n$ a.s. converges to the unique equilibrium point by applying the Kushner-Clark  theorem [c.f. Kushner-Clark, 1978;  Kushner  and Yin, 2003; Duflo, 1997].  But now, it fails to use the Kushner-Clark theorem  because, when $\bm H$ is reducible, for  an equilibrium point $\bm\theta^{\ast}$, its any neighborhood is     not  a region of attraction for $\bm\theta^{\ast}$. However,  the following lemma shows that if a solution  of the ODE
(\ref{eqproofODE}) has  path bounded and bounded from zero, then it must be in  $ \lambda_H \bm S_H$.

\begin{lemma} \label{lem4.0} Suppose Assumption \ref{assumption4.3} is satisfied. Let $\Theta_0=\{\bm\theta\ge \bm 0:\alpha(\bm\theta)\le C_0<0,\; \bm\theta(t) \bm u^t\ge c_0>0\}\supset \lambda_H\bm S_H$. If a solution of  the ODE
(\ref{eqproofODE}) satisfies that $\bm\theta(t)\in \Theta_0$ for all $t\in (-\infty,\infty)$, then $\bm\theta(0)\in \lambda_H \bm S_H$ and $\bm\theta(t)\equiv \bm \theta(0)$.
\end{lemma}

{\bf Proof.} Suppose $\bm\theta(t)$ is a solution of (\ref{eqproofODE}) with the whole path in $\Theta_0$. Let
$f(t)=\int_0^t\frac{1}{\alpha(\bm \theta(s))}ds$. Then $f(t)\to \infty$ and $f(-t)\to -\infty$  as $t\to \infty$, since $\alpha(\bm\theta(s))$ is positive and bounded.
From (\ref{eqproofODE}), it follows that
\begin{equation}\label{eqsolution1} \bm\theta(t)=\bm\theta(0)\exp\left\{-t+f(t)\bm H\right\},
\end{equation}
and
\begin{equation}\label{eqsolution2}
\frac{d}{dt } \bm\theta(t)\bm u^t = -\bm\theta(t)\bm u^t\left(1-\frac{\lambda_H}{\alpha(\bm \theta(t))}\right),
\end{equation}
\begin{equation}\label{eqsolution3} \bm\theta(t)\bm u^t=\bm\theta(0)\bm u^t\exp\left\{-t+\lambda_Hf(t)\right\}.
 \end{equation}

 Let $\bm H$ have the Jordan canonical form (\ref{eqJordanform}), and
\begin{equation}\label{eqsolution4}\bm U=\bm T diag(\bm I_{\nu_1},\bm 0,\cdots,\bm 0)\bm T^{-1}=\sum_{j=1}^{\nu_1}\bm u_j^t \bm v_j, \;\;  \widetilde{\bm H}=\bm H-\lambda_H\bm U.
\end{equation}
Then by (\ref{eqsolution1}) and (\ref{eqsolution3}),
\begin{equation}\label{eqsolution5} \frac{\bm\theta(t)\bm U}{\bm \theta(t)\bm u^t}=\frac{\bm\theta(0)\bm U}{\bm \theta(0)\bm u^t}.
\end{equation}
\begin{equation}\label{eqsolution6} \frac{\bm\theta(t)(\bm I_d-\bm U)}{\bm \theta(t)\bm u^t}=\frac{\bm\theta(0)(\bm I_d-\bm U)}{\bm \theta(0)\bm u^t}\exp\left\{-f(t)(\lambda_H\bm I_d-\widetilde{\bm H})\right\}.
\end{equation}
Note $\bm\theta(t)\in   \Theta$ which implies the right hand of (\ref{eqsolution6}) is bounded,   $f(t)\to -\infty$ as $t\to -\infty$, and that all eigenvalues of $\lambda_H\bm I_d-\widetilde{\bm H}$ have positive real parts. Letting $t\to -\infty$ in (\ref{eqsolution6}) yields $\bm \theta(0)(\bm I_d-\bm U)=\bm 0$, and then $\bm\theta(t)(\bm I_d-\bm U)\equiv \bm 0$. Hence
\begin{equation}\label{eqsolution7} \frac{\bm\theta(t)}{\bm \theta(t)\bm u^t}=\frac{\bm\theta(t)\bm U}{\bm \theta(t)\bm u^t}=\frac{\bm\theta(0)\bm U}{\bm \theta(0)\bm u^t}=\frac{\bm\theta(0)}{\bm \theta(0)\bm u^t},
\end{equation}
\begin{equation}\label{eqsolution8} \frac{\alpha(\bm\theta(t))}{\bm \theta(t)\bm u^t}=\frac{\alpha(\bm\theta(t)\bm U)}{\bm \theta(t)\bm u^t}=\frac{\alpha(\bm\theta(0)\bm U)}{\bm \theta(0)\bm u^t}=\frac{\alpha(\bm\theta(0))}{\bm \theta(0)\bm u^t}\in (0,\infty).
\end{equation}
Combing the above equalities with (\ref{eqsolution2}) yields
$$ \frac{d}{dt} \alpha(\bm\theta(t))=-\alpha(\bm\theta(t)) \left(1-\frac{\lambda_H}{\alpha(\bm \theta(t))}\right)=-\alpha(\bm\theta(t))+\lambda_H. $$
So, $e^t\alpha(\bm\theta(t))-\alpha(\bm\theta(0))=(e^t-1)\lambda_H$. Since $\alpha(\bm\theta(t))$ is bounded, letting $t\to -\infty$ yields $\alpha(\bm\theta(0))=\lambda_H$, and then
$\alpha(\bm\theta(t)) \equiv \lambda_H$.   Now, from (\ref{eqsolution7}) and (\ref{eqsolution8}) it follows that
$$ \bm \theta(t)=\lambda_H\frac{ \bm \theta(t)}{\alpha( \bm \theta(t))}=\lambda_H\frac{\bm\theta(0)\bm U}{\alpha\big(\bm\theta(0)\bm U \big)}\in \lambda_H \bm S_H. \;\;\Box $$

Next, we show that the remainder term $\bm r_{n+1}$ in the stochastic approximation algorithm (\ref{eqSA}) can be neglected. To do so  we need  a  lemma first.
\begin{lemma}\label{lem4.1}
\begin{description}
  \item[\rm (a)] Under Conditions (\ref{th3.5.2}) and (\ref{th3.5.3}) in Theorem \ref{th3.5}, (\ref{eqWLLNM}) is satisfied.
  \item[\rm (b)]   Under  Conditions (\ref{th3.5.1}) and (\ref{eqWLLNM}), for any $T>0$,
\begin{equation}\label{eqlem3.1.2}
\max_{n\le m\le n T}\frac{\alpha(|\bm Y_m|)}{m}\le d\max_{k,q}|H_{k,q}|+o(1),
\end{equation}
\begin{equation}\label{eqlem3.1.3}
\min_{n\le m\le n T}\frac{ \bm Y_m \bm u^t }{m}\ge \lambda_H \min_k u_k+o(1)
\end{equation}
and
\begin{equation}\label{eqlem3.1.4}
\max_{n\le m\le n T}\frac{\alpha(\bm Y_m^-)}{m}\to 0
\end{equation}
in probability.
\item[\rm (c)] Under Conditions (\ref{th3.1.2})-(\ref{th3.1.4}) in Theorem \ref{th3.1}, (\ref{eqSLLNM}) is satisfied.  Under (\ref{th3.1.1}) and (\ref{eqSLLNM}), (\ref{eqlem3.1.2})-(\ref{eqlem3.1.4}) holds almost surely.
  \item[\rm (d)]   Under  Conditions (\ref{th3.3.3}) and (\ref{th3.3.4}) in Theorem \ref{th3.3},  (\ref{eqConvergenceM}) is satisfied.
\end{description}
\end{lemma}

{\bf Proof.} (a)  For (\ref{eqWLLNM}), it is sufficient to show that
\begin{equation}\label{prooflem3.1.1}
\max_{m\le n}\left|\sum_{i=1}^m X_{i,q} (D_{q,k}(i)-H_{q,k}(i))\right|=o(n) \;\; \text{ in probability}.
\end{equation}

By (\ref{th3.5.3}), there is a sequence $0<\epsilon_n\searrow 0$ such that
\begin{equation}\label{prooflem3.1.2}
\frac{1}{\epsilon_n n} \sum_{m=1}^n \ep\left[\|\bm D_m\|I\{\|\bm D_m\|\ge  \epsilon_n  n\}\big|\mathscr{F}_{m-1}\right]\to 0\;\;\text{ in probability}.
\end{equation}
Denote $\xi_i^{(n)}=D_{q,k}(i)I\{\|\bm D_i\|<  \epsilon_n  n\}$. Then
\begin{align*}
& \sum_{i=1}^n  \pr(\xi_i^{(n)}\ne D_{q,k}(i)|\mathscr{F}_{i-1})\le    \sum_{i=1}^n \pr\left(\|D_i\|\ge \epsilon_n n|\mathscr{F}_{i-1}\right) \\
\le & \frac{1}{\epsilon_n n} \sum_{i=1}^n \ep\left[\|\bm D_i\|I\{\|\bm D_i\|\ge  \epsilon_n  n\}\big|\mathscr{F}_{i-1}\right]\to 0\;\;\text{ in probability},
\end{align*}
 by (\ref{prooflem3.1.2}).  So, by Lemma 2.5 of Hall and Heyde (1980),
$$ \pr\left(\xi_i^{(n)}\ne D_{q,k}(i) \text{ for some } i=1,\cdots, n\right)\to 0. $$
Also,
\begin{align*}
&\sum_{i=1}^n |H_{q,k}(n)-\ep[\xi_i^{(n)}|\mathscr{F}_{i-1}]|\\
\le & \sum_{i=1}^n \ep\left[\|\bm D_i\|I\{\|\bm D_i\|\ge  \epsilon_n  n\}\big|\mathscr{F}_{i-1}\right]=o(n)  \;\;\text{ in probability},
 \end{align*}
 by (\ref{prooflem3.1.2}).
Hence, for (\ref{prooflem3.1.1}) it is sufficient to show that
 \begin{equation}\label{prooflem3.1.3}
\max_{m\le n}\left|\sum_{i=1}^m X_{i,q} \big(\xi_i^{(n)}-\ep[\xi_i^{(n)}|\mathscr{F}_{i-1}]\big)\right|=o(n) \;\; \text{ in probability}.
\end{equation}
Note that $\sum_{i=1}^m X_{i,q} \big(\xi_i^{(n)}-\ep[\xi_i^{(n)}|\mathscr{F}_{i-1}]\big)$, $m=1,\cdots, n$, are   martingales with
\begin{align*}
& \sum_{i=1}^n\ep\left[ X_{i,q}^2 \big(\xi_i^{(n)}-\ep[\xi_i^{(n)}|\mathscr{F}_{i-1}]\big)^2|\mathscr{F}_{i-1}\right]\\
\le & \sum_{i=1}^n \epsilon_n \ep\left[|D_{q,k}(i)|\big|\mathscr{F}_{i-1}\right]=\epsilon_n O_P(n)=o(n)\;\;\text{ in probability},
\end{align*}
which implies (\ref{prooflem3.1.3}).  (\ref{eqWLLNM}) is proved.

\medskip
(b)
  By Lemma \ref{lem1} (b),  and the conditions (\ref{th3.5.1}) and (\ref{eqWLLNM}),
 \begin{equation}\label{prooflem3.1.4}
 \max_{n\le m\le n T}\left\|\frac{\bm Y_m}{m}-\frac{1}{m}\sum_{i=1}^m \frac{\bm Y_{i-1}^+}{\alpha(\bm Y_{i-1}^+)} \bm H\right\|\to 0 \;\; \text{ in probability}.
 \end{equation}
 Note that
 $$ \alpha\left(\left|\frac{\bm Y_{i-1}^+}{\alpha(\bm Y_{i-1}^+)} \bm H\right|\right)\le d\max_{q,k}|H_{q,k}| $$
 and
 $$
  \frac{\bm Y_{i-1}^+}{\alpha(\bm Y_{i-1}^+)} \bm H\bm u^t =\frac{\bm Y_{i-1}^+}{\alpha(\bm Y_{i-1}^+)} \lambda_H\bm u^t\ge \lambda_H\min_k u_k. $$
  So, (\ref{eqlem3.1.2}) and (\ref{eqlem3.1.3}) are proved.

  Finally, we show (\ref{eqlem3.1.4}). For $m$ and $k$, let $l_m=\max\{l\le m: Y_{l,k}\ge 0\}$ be the largest integer for which $Y_{l,k}\ge 0$. Then for $n\le m\le n T$,
\begin{align*}
 Y_{m,k}=&Y_{l_m,k}+\sum_{i=l_m+1}^m \sum_{q=1}^dX_{i,q}D_{q,k}(i)\\
 =&Y_{l_m,k}+\sum_{i=l_m+1}^n \sum_{q=1}^dX_{i,q}[D_{q,k}(i)-H_{q,k}(i)]\\
 &+\sum_{i=l_m+1}^m \sum_{q=1}^dX_{i,q}[H_{q,k}(i)-H_{q,k}]
 + \sum_{i=l_m+1}^m \sum_{q=1}^dX_{i,q}H_{q,k}.
 \end{align*}
Note $H_{q,k}\ge 0$ for $q\ne k$. It follows that
\begin{align}\label{eqprooflem3.1.5}
Y_{m,k}^-\le   d\max_{i\le m} \left\|\bm M_{i,2}\right\| +d \sum_{i=1}^m \|\bm H_i-\bm H\|
 + \sum_{i=l_m+1}^m  X_{i,k}|H_{k,k}|.
 \end{align}
For any given $0<\delta<1$, on the event $\big\{\min\limits_{\delta n\le m\le T n}\frac{\bm Y_m\bm u^t}{m}>0\big\}$,  $\min\limits_{\delta n\le m\le T n}\alpha(\bm Y_m^+)>0$. So, if $l_m+1\ge \delta n$, then $X_{i,k}=0$ for $i=l_m+2,\cdots, m$, and if $l_m+1<\delta n$, then $X_{i,k}=0$ for $\delta n+1 \le i\le m$, because $\frac{Y_{i-1,k}^+}{\alpha(\bm Y_{i-1}^+)}=0$ for such $i$. It follows that, on the event $\big\{\min\limits_{\delta n\le m\le T n}\frac{\bm Y_m\bm u^t}{m}>0\big\}$,
 $$ \sum_{i=l_m+1}^m  X_{i,k}|H_{k,k}|\le \delta n |H_{k,k}|. $$
 By  (\ref{eqlem3.1.3}), $\pr\big(\min\limits_{\delta n\le m\le T n}\frac{\bm Y_m\bm u^t}{m}>\lambda_H\min_ku_k/2\big)\to 1$. From (\ref{eqWLLNM}) and (\ref{eqlem3.1.2}) it follows that
$$ \max_{n\le m\le n T} \frac{Y_{m,k}^-}{m}\le \delta +o(1)\;\; \text{in probability}. $$
 The proof of (\ref{eqlem3.1.4}) is completed.

\medskip
(c) By the   strong law of large numbers of martingales,  the condition (\ref{th3.1.4}) implies
$$ \sum_{m=1}^n X_{m,k}\widetilde{D}_{q,k}(m) =o(n) \;\; a.s. $$
The conditions (\ref{th3.1.2}) and (\ref{th3.1.3}) imply
$$ \sum_{m=1}^n X_{m,k}|D_{q,k}(m)-\widetilde{D}_{q,k}(m)|=o(n) \;\; a.s. $$
(\ref{eqSLLNM}) is proved. The reminder of the proof is similar to that of (b).

 \medskip
 (d) It is sufficient to consider each $X_{m,k}(D_{q,k}(m)-H_{q,k}(m))$.    Denote $\overline{D}_{q,k}(m)=D_{q,k}(m)-\widetilde{D}_{q,k}(m)$. Then
\begin{align*}
 &\sum_{m=1}^{\infty} \frac{ \ep\big[|\overline{D}_{q,k}(m)|\big|\mathscr{F}_{m-1}\big]}{m}
  \le     2\frac{ \ep[|D_{q,k}(m)|I\{|D_{q,k}|\ge c_m\}|\mathscr{F}_{m-1}]}{m}<\infty,
 \end{align*}
 by the condition (\ref{th3.3.3}), which implies that
 $\sum_{m=1}^{\infty} \frac{  |\overline{D}_{q,k}(m)|}{m}<\infty \; a.s.$,
 and then  $\sum_{m=1}^{\infty} \frac{  \sum_{i=1}^m|\overline{D}_{q,k}(i)|}{m^2}<\infty \; a.s.$
 For the martingale differences   $\xi_m=:X_{m,k}\widetilde{D}_{q,k}(m)$s,  by the condition (\ref{th3.3.4}) and the stop-time method, without loss of generality, we can assume that
 $$ \sum_{m=1}^{\infty}\frac{ (\log m)^{1+\epsilon}\ep\big[ \widetilde{D}_{q,k}^2(m)\big|\mathscr{F}_{m-1}\big]}{m^2}<M<\infty. $$
 Let $S_m=\sum_{i=1}^m \xi_i$. Then $\ep[\max_{i\le m}|S_i|^2]\le 2\sum_{i=1}^m \ep[ \widetilde{D}_{q,k}^2(m)]$. Hence
 \begin{align*}
 &\ep\left[\sum_{m=1}^{\infty}\frac{\max_{i\le m}|S_i|}{m^2} \right]^2\\
 \le &\ep\left[ \sum_{m=1}^{\infty}\frac{1}{m(\log m)^{1+\epsilon}}
   \sum_{m=1}^{\infty}\frac{(\log m)^{1+\epsilon}\max_{i\le m}|S_i|^2}{m^3} \right]\\
  \le &C \sum_{m=1}^{\infty}\frac{(\log m)^{1+\epsilon}\ep[\max_{i\le m}|S_i|^2]}{m^3} \\
  \le  &C \sum_{m=1}^{\infty}\frac{(\log m)^{1+\epsilon}\sum_{i=1}^m \ep[\widetilde{D}_{q,k}^2(i)]}{m^3}
  \le   C \sum_{i=1}^{\infty} \frac{(\log i)^{1+\epsilon}  \ep[\widetilde{D}_{q,k}^2(i)]}{i^2} \le CM.
 \end{align*}
 Hence
 $$ \sum_{m=1}^{\infty}\frac{\max_{i\le m}|S_i|}{m^2}<\infty\;\; a.s. $$
  The proof is completed. $\Box$

\bigskip
Now, recall that $\bm\theta_n=\frac{\bm Y_n^+}{n}$ satisfies the stochastic approximation algorithm  (\ref{eqSA}) with regression function defined as in (\ref{eqregressF}). For the remainder $\bm r_{m+1}$,   by Lemma \ref{lem4.1}, we have
 for any $T$,
\begin{equation} \label{proofth3.5.4}
\max_{n\le m\le nT}\frac{|\bm s_m|}{m}\to 0,
\end{equation}
 in probability under the assumptions in Theorem \ref{th3.5}, and almost surely under the assumptions in Theorem \ref{th3.1} or \ref{th3.3}, where $\bm s_n=\sum_{m=1}^n \bm r_m$.
  We begin to prove the theorems.
We first prove the weak convergence and then show the almost sure convergence for the non-homogeneity  case. At last, we prove the theorems for the i.i.d. case.

 \bigskip
 {\bf Proof of Theorem \ref{th3.5}.} Note (\ref{proofth3.5.4}), and for any $T$,
\begin{equation} \label{proofth3.5.5} \max_{n\le m\le n T}\alpha(\bm \theta_n)\le d\max_{k,q}|H_{k,q}|+o(1),
\end{equation}
\begin{equation} \label{proofth3.5.6}\min_{n\le m\le n T}  \bm \theta_m \bm u^t  \ge \lambda_H \min_k u_k+o(1)
\end{equation}
 in probability  by (\ref{eqlem3.1.2}) and (\ref{eqlem3.1.3}).

Setting $t_0=0$, $t_n=\frac{1}{1}+\frac{1}{2}+\cdots+\frac{1}{n}$, the 'interpolation function' for a sequence $\{\bm u_n; n\ge 1\}$ in $\mathbb R^d$ is the function from $\mathbb R^+$ to $\mathbb R^d$ defined,  by setting
$$\bm u_0=\bm 0, \;\; \widetilde{\bm u}(t)=\frac{\bm u_n (t_{n+1}-t)+\bm u_{n+1} (t-t_n)}{t_{n+1}-t_n} \;  \text{ for }   t_n\le t\le t_{n+1}, $$
$$\widetilde{\bm u}(t)=\bm u_0 \text{ for }     t\le t_0. $$
 Denote the function which interpolates $\{\bm \theta_n\}$ by $\widetilde{\bm \theta}$ and that which interpolates $\{\bm \sigma_n\}$ where $\bm \sigma_n=\sum_{k=1}^n \frac{\bm r_k}{k}$ by $\widetilde{\bm \sigma}$. We also denote $\overline{\bm \theta}(t)=\bm \theta_n$ if $t_n\le t< t_{n+1}$ and $=\bm u_0$ if $t< t_0$, $h(\bm u_0)=\bm 0$. Then
 $$ \bm\theta(t)=\bm \theta(0)+\int_0^t \bm h\big(\overline{\bm \theta}(s)\big)ds +\widetilde{\bm \sigma}(t), $$
\begin{equation}\label{proofth3.5.7} \widetilde{\bm \theta}(t_n+t)=\widetilde{\bm \theta}(t_n)+\int_0^t \bm h\big(\overline{\bm \theta}(t_n+s)\big)ds +\widetilde{\bm \sigma}(t_n+t)-\widetilde{\bm \sigma}(t_n).
\end{equation}
 We will show that for any $T>0$,
 \begin{equation}\label{proofth3.5.8}
 \max_{t\in [-T,T]}dist\left(\widetilde{\bm \theta}(t_n+t), \lambda_H\bm S_H\right) \to \text{ in probability}.
 \end{equation}
 Let $ m(n,T)=\inf\{k:k\ge n, t_{k+1}-t_n\ge T\}. $ Then $n\le m(n,T)\le 2n e^T$. Note
 \begin{align*}
\sum_{k=n}^j \frac{ \bm r_{k+1}}{k+1}=\sum_{k=n+1}^j \frac{\bm s_k}{k+1}\frac{1}{k}+\frac{\bm s_{j+1}}{j+1}-\frac{\bm s_n}{n+1}.
\end{align*}
It follows that, in probability,
\begin{align*}
&\max_{0\le t\le T}\left\|\widetilde{\bm \sigma}(t_n+t)-\widetilde{\bm \sigma}(t_n)\right\|\le \max_{n\le j\le m(n,T)}\left\|\sum_{k=n}^j \frac{\bm r_{k+1}}{k+1}\right\|\\
\le & \max_{n\le j\le m(n,T)}\frac{\|\bm s_j\|}{j}\sum_{k=n+1}^{m(n,T)} \frac{1}{k}
+2\max_{n\le j\le m(n,T)+1}\frac{\|\bm s_j\|}{j}\\
 \le & (T+3)\max_{n\le j\le 3n e^T}\frac{\|\bm s_j\|}{j}\to 0,
 \end{align*}
 \begin{align*}
 & \max_{0\le t\le T}\alpha\left(\widetilde{\bm \theta}(t_n+t)\right)=\max_{n\le j\le m(n,T)+1}\alpha\left(\bm\theta_j\right)\\
 \le &\max_{n\le j\le 3n e^T}\alpha\left(\bm\theta_j\right)\le d\max_{q,k}|H_{q,k}|+o(1),
 \end{align*}
 $$ \max_{0\le t\le T}  \widetilde{\bm \theta}(t_n+t)\bm u^t\ge \min_{n\le j\le 3n e^T}\bm\theta_j\bm u^t\ge \lambda_H\min_ku_k+o(1)  $$
 and
 \begin{align*}
 & \max_{0\le t\le T} \left\|\widetilde{\bm \theta}(t_n+t)-\overline{\bm\theta}(t_n+t)\right\| \le \max_{n\le j\le m(n,T)} \left\|\bm\theta_{j+1}- \bm\theta_j\right\| \\
 \le   & \max_{n\le j\le 2n e^T } \frac{\|\bm h(\bm \theta_j)\|+\|\bm r_{j+1}\|}{j+1}\le      \max_{n\le j\le 2n e^T } \frac{\|\bm h(\bm \theta_j)\|+\|\bm s_{j+1}\|+\|\bm s_j\|}{j+1} \to 0
 \end{align*}
 by (\ref{proofth3.5.4})-(\ref{proofth3.5.6}). On the other hand, for any give $T$, when $t_n>T$, there is an $n^{\prime}$ and $s_{n^{\prime}}$ with $t_n-T=t_{n^{\prime}}+s_{n^{\prime}}$, $0\le s_{n^{\prime}}\to 0$ and $n^{\prime}\to \infty$. Therefore,  for $t\ge -T$, we have
 $  \widetilde{\bm \theta}(t_n+t)=\widetilde{\bm \theta}(t_{n^{\prime}}+s_{n^{\prime}}+T+t)$,  $ \overline{\bm\theta}(t_n+t) =\overline{\bm \theta}(t_{n^{\prime}}+s_{n^{\prime}}+T+t)$ and
 $$\widetilde{\bm \sigma}(t_n+t)-\widetilde{\bm \sigma}(t_n)=\widetilde{\bm \sigma}(t_{n^{\prime}}+s_{n^{\prime}}+T+t)-\widetilde{\bm \sigma}(t_{n^{\prime}})-\left[\widetilde{\bm \sigma}(t_{n^{\prime}}+s_{n^{\prime}}+T)-\widetilde{\bm \sigma}(t_{n^{\prime}})\right].$$
 It follows that
 \begin{align}\label{proofth3.5.10}
 \max_{t\in[-T, T]}\Big\|\widetilde{\bm \sigma}(t_n+t)-& \widetilde{\bm \sigma}(t_n)\Big\|\to 0,\\
 \label{proofth3.5.11}
 \max_{t\in[-T, T]}\Big\|\widetilde{\bm \theta}(t_n+t)-& \overline{\bm  \theta}(t_n+t)\Big\|\to 0,\\
 \label{proofth3.5.12}
 \max_{t\in[-T, T]}\alpha\left(\widetilde{\bm \theta}(t_n+t)\right) \le & d\max_{q,k}|H_{q,k}|+o(1),
\\
  \label{proofth3.5.13}
    \max_{ t\in[-T, T]}  \widetilde{\bm \theta}(t_n+t)\bm u^t\ge &  \lambda_H\min_ku_k+o(1)
  \end{align}
  in probability for any $T>0$. Write
  \begin{equation}\label{proofth3.5.7ad} \widetilde{\bm \theta}(t_n+t)=\widetilde{\bm \theta}(t_n)+\int_0^t \bm h\big( \widetilde{\bm \theta}(t_n+s)\big)ds +\bm \sigma^{(n)}(t).
\end{equation}
  Then by (\ref{proofth3.5.11}) and the continuity of $\bm h(\cdot)$,
 \begin{equation}\label{proofth3.5.10ad} \max_{t\in[-T, T]}\Big\|\bm\sigma^{(n)}(t)\Big\|\to 0
 \end{equation}
 in probability for any $T>0$.

  Now, we are ready to prove (\ref{proofth3.5.8}). It is sufficient to show for any subsequence $\{n^{\prime}\}$ there is a further subsequence $\{n_j\}\subset \{n^{\prime}\}$ such that (\ref{proofth3.5.8}) holds for this subsequence $\{n_j\}$ almost surely. We can choose the   subsequence $\{n_j\}$ such that (\ref{proofth3.5.10})-(\ref{proofth3.5.13}) and (\ref{proofth3.5.10ad}) holds for this subsequence almost surely for any $T>0$. Except on a null event $\Omega_0$, for each fixed $\omega$, the real functions $\widetilde{\bm \theta}(t_{n_j}+\cdot,\omega)$ and $\bm \sigma^{(n_j)}(\cdot,\omega)$ satisfy (\ref{proofth3.5.12})-(\ref{proofth3.5.10ad}) for any $T>0$.  Hence the sequence of functions $\{\widetilde{\bm \theta}(t_{n_j}+\cdot,\omega)\}$ is bounded and equicontinuous on any bounded interval $[-T,T]$; therefore it relatively compact for the topology of uniform convergence on all compact subsets,  and every limit point $\bm\theta(\cdot)=\bm\theta(\cdot,\omega)$ in sense that
 $$
 \max_{t\in [-T,T]}\left\|\widetilde{\bm\theta}(t_{n_j^{\prime}}+t,\omega)-\bm \theta(t,\omega)\right\| \to 0 \;\; \text{ for any } T>0,
$$
will satisfy the ODE (\ref{eqproofODE})
with the path   $\bm \theta(t)$ $(-\infty<t<\infty)$ in $\Theta=\big\{\bm \theta\ge \bm 0: \alpha(\bm \theta)\le 2 d \max_{q,k}|H_{q,k}|, \bm\theta\bm u^t\ge \lambda_H\min_k u_k/2\big\}$.
By Lemma \ref{lem4.0},   $\bm \theta(t)\equiv\bm\theta(0)\in \lambda_H \bm S_H$. Therefore
     $$
 \max_{t\in [-T,T]}dist\left(\widetilde{\bm\theta}(t_{n_j}+t,\omega), \lambda_H\bm S_H\right) \to 0 \;\; \text{ for any } T>0 \; \text{ and }\omega \not\in \Omega_0.
 $$
 That is
   $$
 \max_{t\in [-T,T]}dist\left(\widetilde{\bm\theta}(t_{n_j}+t), \lambda_H\bm S_H\right) \to 0 \;\; a.s. \;\; \text{ for any } T>0,
 $$
and (\ref{proofth3.5.8}) is proved.

Now,  note  for $T>1$,
 $$\max\limits_{n\le m\le nT}dist\Big(\bm\theta_m, \lambda_H\bm S_H\Big)\le \max\limits_{t\in [0,T]}dist\Big(\bm\theta(t_n+t), \lambda_H\bm S_H\Big). $$
(\ref{th3.5.4}) is proved. (\ref{th3.5.5}) follows from (\ref{th3.5.4}) immediately. Finally, note
$$ \frac{\bm N_n}{n}=\frac{1}{n} \sum_{m=1}^n \frac{\bm Y_{m-1}^+}{\alpha(\bm Y_{m-1}^+)}+\frac{\bm M_{n,1}}{n}. $$
So, for $n\le m\le n T$,
\begin{align*}
& dist\Big(\frac{\bm N_m}{m}, \bm S_H\Big)\le \frac{1}{m}\sum_{i=1}^n dist\Big(\frac{\bm Y_{i-1}^+}{\alpha(\bm Y_{i-1}^+)}, \bm S_H\Big)+\frac{\|\bm M_{m,1}\|}{m}\\
\le & \frac{1}{m}\sum_{i=1}^{\delta n} 2 + \frac{1}{m}\sum_{i=\delta n+1}^n \max_{\delta n\le i\le n T }dist\Big(\frac{\bm Y_i^+}{\alpha(\bm Y_i^+)}, \bm S_H\Big)+\frac{\|\bm M_{m,1}\|}{m} \\
\le & 2\delta + \max_{\delta n\le i\le n T }dist\Big(\frac{\bm Y_i^+}{\alpha(\bm Y_i^+)}, \bm S_H\Big)+ \max_{ n\le m\le n T }\frac{\|\bm M_{m,1}\|}{m}.
\end{align*}
(\ref{th3.5.6}) follows from (\ref{th3.5.5}) immediately.
$\Box$

\bigskip
{\bf Proof of Theorem \ref{th3.1}.}   Now,  by Lemma \ref{lem4.1} (c) (\ref{proofth3.5.4})-(\ref{proofth3.5.10ad}) holds almost surely for any $T>0$. So, almost surely, every limit point $\bm\theta(\cdot)$ of the sequence $\{\widetilde{\bm\theta}(t_n+\cdot)\}$ will satisfy $\bm\theta(t)\in \lambda_H\bm S_H$, which implies
$$ \lim_{n\to \infty} dist\left(\frac{\bm Y_n^+}{n},\lambda_H\bm S_H\right)=\lim_{n\to \infty} dist\left(\widetilde{\bm\theta}(t_n),\lambda_H\bm S_H\right)=0\;\; a.s. $$
The proof is completed by noting $\frac{\bm Y_n^-}{n}\to 0$ a.s.
$\Box$

\bigskip
{\bf Proof of Theorem \ref{th3.3}.} (\ref{eqSLLNM}) remains true. (\ref{th3.3.1}) implies (\ref{th3.1.1}) by (\ref{eqConvergenceM}).  By Theorem \ref{th3.1}, (\ref{eqth2.1.1}) holds. Let $\bm U=\sum_j\bm u_j^t\bm v_j$ be defined as in (\ref{eqsolution4}), then $\frac{\bm\theta\bm U}{\alpha(\bm \theta\bm U)}$ is a project from $\mathbb R^d_+$ to $\bm S_H$. So, by (\ref{eqth2.1.1}),
$$ \bm\theta_n-\lambda_H\frac{\bm\theta_n\bm U}{\alpha(\bm \theta_n\bm U)}\to 0 \;\; a.s. $$
It is sufficient to show that
$$ \frac{\bm\theta_n\bm U}{\alpha(\bm \theta_n\bm U)} \text{ converges } a.s. $$
By the stochastic approximation algorithm (\ref{eqSA}),
\begin{equation}\label{eqproofth3.3.0}\bm\theta_{n+1}\bm U-\bm\theta_n\bm U=-\frac{\left(1-\lambda_H/\alpha(\bm \theta_n)\right)\bm \theta_n\bm U}{n+1} +\frac{\bm r_{n+1}\bm U}{n+1}.
\end{equation}
Write $\bm x_n=\bm\theta_n\bm U$ and $\bm\delta_n=\bm s_n\bm U/n$. We can rewrite the above equality as
$$ \bm x_{n+1}-\bm\delta_{n+1} -\left(\bm x_n- \bm \delta_n \right)=-\frac{\left(1-\lambda_H/\alpha(\bm \theta_n)\right)\bm x_n}{n+1} +\frac{\bm \delta_n}{n+1}. $$
 It is obvious that $\bm\delta_n\to \bm 0$  a.s. Let $f(\bm x)=\frac{\bm x}{\alpha(\bm x)}$. Then $\frac{\partial f(\bm x)}{\partial \bm x}=\frac{\bm I_d}{\alpha(\bm x)}-\frac{\bm 1^t\bm x}{\alpha^2(\bm x)}$, $\bm x\frac{\partial f(\bm x)}{\partial \bm x}\equiv \bm 0$. Note that $f(\bm x)$, $\frac{\partial f(\bm x)}{\partial \bm x}$ and the derivations of $\frac{\partial f(\bm x)}{\partial \bm x}$ are all bounded on a neighborhood of $\bm S_H$.
 Let $\bm x_n^{\prime}=\bm x_n-\bm\delta_n$. It follows that
\begin{align}\label{eqproofth3.3.1}
&f\left(\bm x_{n+1}-\bm\delta_{n+1}  \right)-f\left(\bm x_n-\bm\delta_n\right)
\nonumber\\
=&\left(-\frac{\left(1-\lambda_H/\alpha(\bm \theta_n)\right)\bm x_n}{n+1} +\frac{\bm \delta_n}{n+1}\right) \frac{\partial f(\bm x)}{\partial \bm x}\Big|_{\bm x=\bm x_n-\bm \delta_n}
\nonumber\\
 &+  \frac{O(1)}{(n+1)^2}+ \frac{O(\|\bm\delta_n\|^2)}{(n+1)^2}
 \nonumber  \\
 =& -\frac{\left(1-\lambda_H/\alpha(\bm \theta_n)\right)}{n+1}  \bm x_n\frac{\partial f(\bm x)}{\partial \bm x}\Big|_{\bm x=\bm x_n}
 +\frac{O(1)}{(n+1)^2}+ \frac{O(\|\bm\delta_n\|)}{n+1}
 \nonumber \\
 &-\frac{\left(1-\lambda_H/\alpha(\bm \theta_n)\right)}{n+1}   \bm x_n\left(\frac{\partial f(\bm x)}{\partial \bm x}\Big|_{\bm x=\bm x_n-\bm \delta_n}-\frac{\partial f(\bm x)}{\partial \bm x}\Big|_{\bm x=\bm x_n}\right)
 \nonumber\\
 =&  \frac{O(1)}{(n+1)^2}+ \frac{O(\|\bm\delta_n\|)}{n+1}.
\end{align}
 We will prove that
\begin{equation} \label{eqdelta1}
  \sum_{m=1}^{\infty} \frac{\|\bm \delta_m\|}{m+1}<\infty\;\; a.s.
 \end{equation}
Then $f(\bm x_n-\bm\delta_n)$ is convergent a.s., which implies that $f(\bm x_n)$ is convergent a.s. because $\bm\delta_n\to \bm 0$ a.s.

 Note $\liminf_{n\to \infty}\bm Y_n^+\bm u^t/n>0$ a.s.  Without loss of generality, we can assume that $\alpha(\bm Y_i^+)\ge c\bm Y_i^+\bm u^t\ne 0$ a.s. for all $i$. By (\ref{eqprooflem3.1.5}),
$$Y_{m,k}^-\le \max_{i\le m}\|\bm M_{i,2}\|+\sum_{i=1}^m\|\bm H_i-\bm H\|+X_{l_m+1,k}|H_{k,k}|. $$
So,
$$\|\bm s_m\|\le \max_{i\le m}\|\bm M_{i,2}\| +C\max_{i\le m}\|\bm M_{i,1}\|  +C\sum_{i=1}^m\|\bm H_i-\bm H\|+C\; a.s.
$$
Note
$$\sum_{m=1}^{\infty} \frac{\sum_{i=1}^m\|\bm H_i-\bm H\|}{m(m+1)}=\sum_{i=1}^{\infty}\sum_{m=1}^{\infty} \frac{\|\bm H_i-\bm H\|}{m(m+1)}\le \sum_{i=1}^{\infty} \frac{|\bm H_i-\bm H\|}{i}<\infty\;\;a.s. $$
by the assumption (\ref{th3.3.1}),
$$\sum_{m=1}^{\infty} \frac{\max\limits_{i\le m}\|\bm M_{i,1}\|}{m(m+1)}\le C\sum_{m=1}^{\infty} \frac{\sqrt{m\log\log m}}{m(m+1)} <\infty \;\; a.s. $$
by Lemma \ref{lem1} (b),
and the fact (\ref{eqConvergenceM}).  So, (\ref{eqdelta1}) is satisfied. $\Box$

\bigskip
{\bf Proof of Theorem \ref{th3.4}.}
 Under the assumptions in the theorem,    $\bm \theta_n=\frac{\bm Y_n^+}{n}$ satisfies the stochastic approximation algorithm  (\ref{eqSA}) with
\begin{align*} \|\bm s_n\|\le & C\max_{m\le n}\|\bm M_{n,1}\|+C\max_{m\le n}\|\bm M_{n,2}\|+C\sum_{m=1}^n\|\bm H_m-\bm H\|+C\\
=& O\big(\sqrt{n\log\log n}\big)\;\; a.s.
\end{align*}
By (\ref{eqproofth3.3.1}),
\begin{align*}
&\|f(\bm x_n-\bm\delta_n)-\bm V\|\le   \sum_{m=n-1}^{\infty}\left(\frac{O(1)}{(m+1)^2}+\frac{O(\|\bm\delta_m\|)}{m(m+1)} \right)\\
\le & \sum_{m=n-1}^{\infty}\frac{O\big(\sqrt{m\log\log m}\big)}{m(m+1)}=O\Big(n^{-1/2}\sqrt{\log\log n}\Big)\;\; a.s.,
\end{align*}
where $\bm x_n=\bm\theta_n\bm U$, $\bm\delta_n=\bm s_n\bm U/n$ and $f(\bm x)=\frac{\bm x}{\alpha(\bm x)}$.
It follows that
\begin{align}\label{eqproofth2.3.3} \|f(\bm x_n)-\bm V\|\le & \|f(\bm x_n-\bm \delta_n)-\bm V\|+O(\|\bm\delta_n\|)\nonumber\\
=& O\Big(n^{-1/2}\sqrt{\log\log n}\Big)\;\; a.s.
\end{align}
By the stochastic approximation algorithm (\ref{eqSA}),
$$ \bm\theta_{n+1}(\bm I_d-\bm U)=\bm \theta_n(\bm I_d-\bm U)\left(\bm I_d-\frac{\bm I_d- \widetilde{\bm H}/\alpha(\bm \theta_n)}{n+1}\right)+\frac{\bm r_{n+1}(\bm I-\bm U)}{n+1}. $$
Let $\bm H_{j+1}=\bm I_d-\widetilde{\bm H}/\alpha(\bm\theta_i)$ and
\begin{equation}  \bm \Pi_m^n  =\Big(\bm I_d-\frac{\bm H_{m+1}}{m+1}\Big)\cdots \Big(\bm I_d-\frac{\bm H_{n}}{n}\Big),\; \widetilde{\bm \Pi}_m^n=\prod_{j=m+1}^n\Big(\bm I_d-\frac{\bm I_d-\frac{\widetilde{\bm H}}{\lambda_H}}{j}\Big).
\end{equation}
  Then $\bm H_j\to \bm I_d-\widetilde{\bm H}/\lambda_H$ a.s. as $j\to \infty$ because we have shown that $\alpha(\bm \theta_j)\to\lambda_H\alpha(\bm V)=\lambda_H$ a.s. in Theorem \ref{th3.3}. Note that the smallest real part of the eigenvalues of  $\bm I_d-\widetilde{\bm H}/\lambda_H$ is $1-\rho>0$. It follows that  $\|\bm \Pi_m^n\|\le C_{\delta}\big(n/m)^{\rho-1+\delta}$   for any $\delta>0$  and $\|\widetilde{\bm \Pi}_m^n\|\le C_{\delta}\big(n/m)^{\rho-1}(\log n/m)^{\nu_{sec}-1}$ by Lemma B.1   of Zhang (2016). Similar to (2.17) and (2.22) of Zhang (2016),
we have for $\delta>0$ small enough,
 \begin{align}\label{eqproofth2.3.4}
 & \bm \theta_n-\bm x_n=\bm \theta_n(\bm I_d-\bm U) =  \bm \theta_0(\bm I_d-\bm U) \bm \Pi_0^n+ \sum_{m=1}^n \frac{\bm r_m(\bm I_d-\bm U)}{m}\bm \Pi_m^n \nonumber\\
=&\bm \theta_0(\bm I_d-\bm U) \bm \Pi_0^n+\frac{\bm s_n(\bm I_d-\bm U)}{n}\bm \Pi_n^n+\sum_{m=1}^{n-1}\bm s_m(\bm I_d-\bm U)\frac{\bm I_d-\bm H_{m+1}}{m(m+1)}\bm \Pi_{m+1}^n\nonumber\\
=&O\big(n^{\rho-1+\delta}\big)  +\frac{O(\sqrt{n\log\log n})}{n} +\sum_{m=1}^{n-1}\frac{O(\sqrt{m\log\log m})}{m(m+1)}\big(\frac{n}{m}\big)^{\rho-1+\delta}\nonumber\\
=&\begin{cases} O(n^{\rho-1+\delta}) \;\; a.s. & \text{ if } \rho\ge 1/2, \\
O\Big(n^{-1/2}\sqrt{\log\log n}\Big)\;\; a.s. & \text{ if } \rho< 1/2
\end{cases}
=:O(\overline{b}_n) \;\; a.s.
\end{align}
So,
\begin{equation}\label{eqproofth2.3.6} f(\bm \theta_n)-f(\bm x_n)=O(\|\bm\theta_n-\bm x_n\|)=O(\overline{b}_n)\;\; a.s.
\end{equation}
Combining (\ref{eqproofth2.3.3}) and (\ref{eqproofth2.3.6}) yields
$$ \frac{\bm \theta_n}{\alpha(\bm\theta_n)}-\bm V =f(\bm \theta_n)-f(\bm x_n)+(f(\bm x_n)-\bm V)=O(\overline{b}_n)\;\; a.s. $$
By the stochastic approximation algorithm  (\ref{eqSA}) again,
\begin{align*}
 &\bm\theta_n-\lambda_H\bm V= \frac{1}{n}\sum_{m=1}^n\left(\frac{\bm \theta_{m-1}}{\alpha(\bm\theta_{m-1})}-\bm V\right)\bm H+\frac{\bm s_n}{n} \\
 =&\frac{1}{n}\sum_{m=1}^nO(\overline{b}_m) +O\big(n^{-1/2}\sqrt{\log\log n}\big) =O(\overline{b}_n)\;\; a.s.
 \end{align*}
Hence
$$ \alpha(\bm\theta_n)-\lambda_H=O(\overline{b}_n)\;\; a.s.$$
Then $H_{j+1}=\bm I_d-\widetilde{\bm H}/\lambda_H+O(\overline{b}_j)$ a.s., and
  we can rewrite the stochastic approximation algorithm  for $\bm \theta_n(\bm I-\bm U)$ as
$$ \bm\theta_{n+1}(\bm I_d-\bm U)=\bm \theta_n(\bm I_d-\bm U)\left(\bm I_d-\frac{\bm I_d- \widetilde{\bm H}/\lambda_H}{n+1}\right)+\frac{\bm r_{n+1}^*}{n+1}, $$
with $\bm r_{n+1}^*=\bm r_{n+1}(\bm I_d-\bm U)+O(\overline{b}_n )$ a.s.  Repeating the above arguments  yields
\begin{align*}
   \bm \theta_n-\bm x_n
=& \bm \theta_0(\bm I_d-\bm U)  \widetilde{\bm \Pi}_0^n+\frac{\bm s_n^*}{n}\widetilde{\bm \Pi}_n^n+\sum_{m=1}^{n-1}\bm s_m^* \frac{\widetilde{\bm H}/\lambda_H}{m(m+1)}\widetilde{\bm \Pi}_{m+1}^n\nonumber\\
=&O\big(n^{\rho-1}(\log n)^{\nu_{sec}-1}\big)  +\frac{O(n\overline{b}_n^2+\sqrt{n\log\log n})}{n} \nonumber\\
& +\sum_{m=1}^{n-1}\frac{O(m\overline{b}_m^2+\sqrt{m\log\log m})}{m(m+1)}\big(\frac{n}{m}\big)^{\rho-1}(\log n/m)^{\nu_{sec}-1}\nonumber\\
=&\begin{cases} O\big(n^{\rho-1}(\log n)^{\nu_{sec}-1}\big) \;\; a.s. & \text{ if } \rho> 1/2, \\
O\big(n^{-1/2}(\log n)^{\nu_{sec}-1}(\log\log n)^{1/2}\big) \;\; a.s. & \text{ if } \rho= 1/2, \\
O\Big(n^{-1/2}\sqrt{\log\log n}\Big)\;\; a.s. & \text{ if } \rho< 1/2
\end{cases}\\
=&O(b_n) \;\; a.s.,
\end{align*}
$$ \frac{\bm \theta_n}{\alpha(\bm\theta_n)}-\bm V =f(\bm \theta_n)-f(\bm x_n)+(f(\bm x_n)-\bm V)=O(b_n)\;\; a.s. $$
and
$$ \bm\theta_n-\lambda_H\bm V =O(b_n)\;\; a.s.$$
Note $\|\bm Y_n^-\|=O(\sqrt{n\log\log n})$ a.s. and
$$ \frac{\bm N_n}{n}-\bm V= \frac{1}{n}\sum_{m=1}^n\left(\frac{\bm \theta_{m-1}}{\alpha(\bm\theta_{m-1})}-\bm V\right) +\bm M_{n,1}. $$
The proof is completed. $\Box$

\bigskip
The theorems for non-homogeneity case have been proved. Finally, we show the main results in Section \ref{sectionmain}.

{\bf Proof of Theorem \ref{th2.1}.} Note that $\{\bm D_n\}$ are i.i.d.  and (\ref{th3.1.1}) is obvious   since $\bm H_n=\bm H$. Also,    $\ep[|D_{k,q}(m)|]<\infty$   implies (\ref{th3.1.2})-(\ref{th3.1.4}) with $c_m=m$.      Theorem \ref{th2.1} follows from Theorem \ref{th3.1} immediately.

\bigskip
{\bf Proof of Theorem \ref{th2.3}.} In Remark \ref{remark3.2}, we have shown that $\ep\big[\|\bm D_n\|\log(\|\bm D_n\|)\big]<\infty$ implies (\ref{th3.3.3}) and (\ref{th3.3.4})  in Theorem \ref{th3.3} with $c_m=m/(\log m)^{1+\epsilon}$. So, (\ref{eqth2.3.1})-(\ref{eqth2.3.3}) hold.

Now, assume $D_{k,q}(m)\ge 0$ for all $q,k$ and $n$. The goal is to   show that $\varpi_j>0$. Note that $\bm u_j^t$ is a  right eigenvector of $\bm H$ corresponding to $\lambda_H$. Then
$$ \ep[\bm Y_{n+1}\bm u_j^t|\mathscr{F}_n]]=\bm Y_n\bm u_j^t\Big(1+\frac{\lambda_H}{\alpha(\bm Y_n)}\Big)\le \bm Y_n\bm u_j^t\exp\Big\{ \frac{\lambda_H}{\alpha(\bm Y_n)}\Big\}. $$
It follows that $\bm Y_n\bm u_j^t\exp\{-\lambda_Hq_{n-1}\}$ is a nonnegative supermartingale, where $q_n= \sum_{l=0}^{n-1}\frac{1}{\alpha(\bm Y_l)}$, and so it  converges to a nonnegative  random variable, say $\varphi_j$, almost surely. On the other hand, note
$$ \frac{\bm Y_n\bm u_j^t}{\bm Y_{n+1}\bm u_j^t} =1 -\frac{\bm X_{n+1}\bm D_{n+1}\bm u^t_j}{\bm Y_n\bm u_j^t}+\frac{(\bm X_{n+1}\bm D_{n+1}\bm u^t_j)^2}{\bm Y_n\bm u_j^t(\bm Y_n\bm u_j^t+\bm X_{n+1}\bm D_{n+1}\bm u^t_j)}.$$
It follows that
$$
\ep\left[\frac{\bm Y_n\bm u_j^t}{\bm Y_{n+1}\bm u_j^t}\Big|\mathscr{F}_n\right]=1-\frac{\lambda_H}{\alpha(\bm Y_n)}+\frac{1}{\bm Y_n\bm u_j^t}\sum_{k=1}^d \frac{Y_{n,k}}{\alpha(\bm Y_n)}
\ep\left[\frac{(\bm d_1^{(k)}\bm u_j^t)^2}{\bm Y_n\bm u_j^t+\bm d_1^{(k)}\bm u_j^t}\Big|\mathscr{F}_n\right].
$$
where $\bm d_1^{(k)}$ is the $k$-th row of $\bm D_n$.
Since $\ep[\bm D_n\bm u_j^t]=\lambda_H\bm u_j^t$, if the $k$-th element $u_j^{(k)}$ of $\bm u_j$ is zero, then $\ep[\bm d_1^{(k)}\bm u_j^t]=u_j^{(k)}=0$, and then $\bm d_1^{(k)}\bm u_j^t$ must be zero because it is a nonnegative random variable. So, the summation over $k$ is taken over those $k$s for which the $k$-th element of $\bm u_j$ is not zero. Then $\sum_k Y_{n,k}\le c \bm Y_n\bm u_j^t$ where $c=\max\{1/\mu_j^{(k)}: \mu_j^{(k)}\ne 0, k=1,\cdots, d\}$. It follows that
\begin{align*}
\ep\left[\frac{\bm Y_n\bm u_j^t}{\bm Y_{n+1}\bm u_j^t}\Big|\mathscr{F}_n\right]
\le & 1-\frac{\lambda_H}{\alpha(\bm Y_n)}+c\frac{1}{\alpha(\bm Y_n)}\max_k
\ep\left[\frac{(\bm d_1^{(k)}\bm u_j^t)^2}{\bm Y_n\bm u_j^t+\bm d_1^{(k)}\bm u_j^t}\Big|\mathscr{F}_n\right]
\\
\le & \exp\left\{- \frac{\lambda_H}{\alpha(\bm Y_n)}+c\frac{1}{\alpha(\bm Y_n)}\max_kf_k(\bm Y_n\bm u_j^t, \bm u_j)\right\},
\end{align*}
 where $f_k(l,\bm y)=\ep\left[\frac{(\bm d_1^{(k)}\bm y^t)^2}{l+\bm d_1^{(k)}\bm y^t}\right]$.    It follows that
$$ \frac{1}{\bm Y_n\bm u_j^t}\exp\left\{ \lambda_H q_{n-1}-\sum_{m=0}^{n-1} c\max_k\frac{f_k(\bm Y_m\bm u_j^t, \bm u_j)}{\alpha(\bm Y_m)}\right\} $$
 is a nonnegative supermartingale, and so it   also converges to a nonnegative  random variable   almost surely. Next, we will show that
$\sum_{m=0}^{\infty} \frac{f_k(\bm Y_m\bm u_j^t, \bm u_j)}{\alpha(\bm Y_m)}<\infty$ a.s. Then both $\frac{1}{\bm Y_n\bm u^t}\exp\left\{ \lambda_H q_{n-1}\right\}$ and $ \bm Y_n\bm u^t \exp\left\{- \lambda_H q_{n-1}\right\}$ converge to  nonnegative random variables.
Hence $\pr(\varphi_j>0)=1$.

 First, note
 $$ \sum_{m=2}^{\infty}\pr(X_{m+1,k}=1|\mathscr{F}_m)\ge \sum_{m=2}^{\infty} \frac{Y_{0,k}}{\alpha(\bm Y_{m-1})}\ge c \sum_{m=1}^{\infty}\frac{1}{m}=+\infty. $$
So, $\pr(X_{m,k}=1\; i.o.)=1$, which implies $N_{n,k}=\sum_{m=1}^nX_{m,k}\to \infty$ a.s. If the $k$-th element $u_j^{(k)}$  of $\bm u_j$ is not zero, then $\ep[\bm d_1^{(k)}\bm u_j^t]=\lambda_H u_j^{(k)}>0$. Hence,   $\bm Y_n\bm u_j^t\ge \sum_{m=1}^n X_{m,k}\bm d_m^{(k)}\bm u_j^t\sim N_{n,k}\ep[\bm d_1^{(k)}\bm u_j^t]\to \infty$ a.s.

Next, from $\bm Y_n\bm u_j^t\to \infty$ a.s., we conclude that $f_k(\bm Y_n\bm u_j^t,\bm u_j)\to 0$ a.s. and  then
$$ c\exp\{ -\lambda_H(1-\delta)q_{n-1}\}\le \bm Y_n\bm u_j^t\le C\exp\{-\lambda_H q_{n-1}\}. $$
Taking the summation over $j$ and noting that $(\sum_{j=1}^{\nu_1} \bm Y_n\bm u_j^t)/n\to \lambda_H$ a.s., yields that
$   c\exp\{ -\lambda_H(1-\delta)q_{n-1}\}\le n \le C\exp\{-\lambda_H q_{n-1}\}$ a.s. and
$$ cn^{1-\delta}\le \bm Y_n\bm u_j^t \le Cn \; a.s. $$

 Finally, by noting
 $ f_k(\bm Y_n\bm u_j^t,\bm u_j)\le \frac{n^{1/2}}{\bm Y_n\bm u_j^t}+\max_k \ep[  \bm d_1^{(k)}\bm u_j^tI\{\bm d_1^{(k)}\bm u_j^t\ge n^{1/2}\}]$ and $\frac{\alpha(\bm Y_n)}{n}\to \lambda_H$ a.s., we have
 \begin{align*}
 &\sum_m \frac{1}{\alpha(\bm Y_m)}f(\bm Y_n\bm u_j^t,\bm u_j)\\
 \le & C\sum_m \frac{m^{1/2}}{m m^{1-\delta}}+C\sum_m \frac{\max_k \ep[  \bm d_1^{(k)}\bm u_j^tI\{\bm d_1^{(k)}\bm u_j^t\ge m^{1/2}\}]}{m}\\
 \le & C+C\max_k\ep[(\bm d_1^{(k)}\bm u_j)\log (\bm d_1^{(k)}\bm u_j)]<\infty.
 \end{align*}
 So, $\bm Y_n\bm u_j^t\exp\{-\lambda_H q_{n-1}\}\to \varphi_j>0$ a.s., which implies that
 $$ \frac{\bm Y_n\bm u_j^t/n}{\sum_{i=1}^{\nu_1} \bm Y_n\bm u_i^t/n}=\frac{\bm Y_n\bm u_j^t}{\sum_{i=1}^{\nu_1} \bm Y_n\bm u_i^t}\to \frac{\varphi_j}{\sum_{i=1}^{\nu_1} \varphi_i}\;\; a.s.  $$
 We conclude that
 $$\varpi_j= \bm V\bm u_j^t  =\frac{\varphi_j}{\sum_{i=1}^{\nu_1} \varphi_i}>0\;\; a.s., \;\; j=1,\cdots, \nu_1. $$

 For showing that $\varphi_j$ has no point probability mass in $[0,1)$ we apply the conditional central limit theorem. Now we have a more condition that $\ep[\|\bm D_n\|^2]<\infty$.  Consider (\ref{eqproofth3.3.0}) again, where $\bm r_{n+1}=\Delta \bm M_{n+1,1}\bm H+\Delta \bm M_{n+1,2}$ since $\bm H_n=\bm H$ and $\bm Y_n\ge \bm Y_0>\bm 0$.
 Similar to (\ref{eqproofth3.3.1}) for $\bm x_n=\bm\theta_n\bm U$ and $f(\bm x)=\bm x/\alpha(\bm x)$,  we have
 \begin{align*}
 f(\bm x_{n+1})-f(\bm x_n)=&\left(-\frac{\left(1-\lambda_H/\alpha(\bm \theta_n)\right)\bm x_n}{n+1} +\frac{\bm r_{n+1}\bm U}{n+1}\right) \frac{\partial f(\bm x)}{\partial \bm x}\Big|_{\bm x=\bm x_n}
 \\
 &+  \frac{O(1)}{(n+1)^2}+ \frac{O(\|\bm r_{n+1}\|^2)}{(n+1)^2}
 \\
 =& \frac{\bm r_{n+1}\bm U}{n+1}  \frac{\partial f(\bm x)}{\partial \bm x}\Big|_{\bm x=\bm x_n}
  +  \frac{O(1)}{(n+1)^2}+ \frac{O(\|\bm D_{n+1}\|^2)}{(n+1)^2}.
 \end{align*}
 It is easily shown that $\sqrt{n} \sum_{m=n}^{\infty}\frac{\|\bm D_{m}\|^2}{m^2}=o(1)$ a.s. Write $\bm \eta_n = \bm r_n\bm U\frac{\partial f(\bm x)}{\partial \bm x}\Big|_{\bm x=\bm x_{n-1}}$. Then $\{\bm \eta_n\}$ is a sequence of martingale differences with
 \begin{align*}
  & \ep[\bm \eta_{n+1}^t    \bm \eta_{n+1}|\mathscr{F}_n]\\
  =& \Big(\bm U\frac{\partial f(\bm x)}{\partial \bm x}\Big|_{\bm x=\bm x_{n}}\big)^t
 \left(\bm H^t\Big(diag\big(\frac{\bm \theta_n}{\alpha(\bm\theta_n)}\big)-\frac{\bm \theta_n^t}{\alpha(\bm\theta_n)}   \frac{\bm \theta_n}{\alpha(\bm\theta_n)} \Big)\bm H\right.\\
 & \;\;  \left.+\ep_{\bm D_1}\big[(\bm D_1-\bm H)^t diag \big(\frac{\bm \theta_n}{\alpha(\bm\theta_n)}\big)(\bm D_1-\bm H) \right)\bm U\frac{\partial f(\bm x)}{\partial \bm x}\Big|_{\bm x=\bm x_{n}}\\
 \to & \Big(\bm U (\bm I_d-\bm 1^t\bm V)^t/\lambda_H
 \left[\bm H^t(diag(\bm V)-\bm V^t\bm V)\bm H\right.\\
 & \;\; \left.+\ep_{\bm D_1}\big[(\bm D_1-\bm H)^t diag (\bm V)(\bm D_1-\bm H)\big] \right]\bm U (\bm I_d-\bm 1^t\bm V)/\lambda_H\\
 = &   (\bm I_d - \bm V^t\bm 1 )^t\bm U^t
  diag(\bm V)\bm U (\bm I_d -\bm 1^t\bm V)  \\
 & +  (\bm I_d - \bm V^t \bm 1) \bm U^t\ep_{\bm D_1}\big[(\bm D_1-\bm H)^t diag (\bm V)(\bm D_1-\bm H)\big]   \bm U (\bm I_d-\bm 1^t\bm V)/\lambda_H^2\\
  =: & \bm\Sigma_1+\bm\Sigma_2=\bm\Sigma \;\; a.s.
 \end{align*}
  It follows that
  $$n\sum_{m=n}^{\infty}\ep\left[\frac{\bm\eta_m^t\bm\eta_m}{m^2}|\mathscr{F}_{m-1}\right]\to \bm\Sigma\; a.s. $$
  Further, the conditional Lindeberger condition is satisfied, i.e., for any $\epsilon>0$,
 \begin{align*}
 &\sum_{m=n}^{\infty}\ep\left[\frac{n\|\bm\eta_m\|^2}{m^2}I\left\{\frac{\sqrt{n}\|\bm\eta_m\|^2}{m}\ge \epsilon\right\} |\mathscr{F}_{m-1}\right]\\
 \le & c\sum_{m=n}^{\infty} \frac{n }{m^2}\ep\left[\|\bm D_1\|^2I\{\|\bm D_1\|\ge \epsilon \sqrt{n}\}\right]\ \to 0\;\; a.s.
 \end{align*}
  With the above results,  we can show that given $\mathscr{F}_n$, the conditional distribution of $\sqrt{n}\sum_{m=1}^{\infty}\frac{\bm\eta_m}{m}$ will almost surely converge to a multi-normal distribution
  $N(\bm 0, \bm\Sigma(\omega))$. It follows that
  $$ \ep\left[e^{it\sqrt{n}(\bm V-f(\bm x_n))\bm u_j^t}\big|\mathscr{F}_n\right]\to e^{-\frac{t^2}{2}\sigma_j^2} \;\; a.s., $$
  where $\sigma_j^2=\bm u_j\bm\Sigma\bm u_j^t$.
 For any $0<p<1$, let $E=\{\varphi_j=p\}=\{\bm V\bm u_j^t=p\}$, $I_n=\ep[I_E|\mathscr{F}_n]$. Then $I_n\to I_E$ a.s., $\ep[|I_n-I_E|\big|\mathscr{F}_n]\to 0$ in $L_1$, and so
  \begin{align*}
  &\lim_n\ep\left[e^{it\sqrt{n}(\bm V-f(\bm x_n))\bm u_j^t}I_E\big|\mathscr{F}_n\right]
  =\lim_n\ep\left[e^{it\sqrt{n}(\bm V-f(\bm x_n))\bm u_j^t}I_n\big|\mathscr{F}_n\right]\\
  =&\lim_n\ep\left[e^{it\sqrt{n}(\bm V-f(\bm x_n))\bm u_j^t}\big|\mathscr{F}_n\right]I_n=e^{-\frac{t^2}{2}\sigma_j^2} I_E\;  \text{ in } L_1.
  \end{align*}
  Note on the event $E$, $(\bm V-f(\bm x_n))\bm u_j^t=p-f(\bm x_n)\bm u_j^t$ is $\mathscr{F}_n$ measurable. We conclude
  \begin{align*}
  I_E=\lim_n \ep[I_E|\mathscr{F}_n]= \lim_n\left|\ep\left[e^{it\sqrt{n}(\bm V-f(\bm x_n))\bm u_j^t}I_E\big|\mathscr{F}_n\right]\right|=e^{-\frac{t^2}{2}\sigma_j^2} I_E\; \text{ in } L_1.
  \end{align*}
  So, $I_E=e^{-\frac{t^2}{2}\sigma_j^2} I_E$ a.s. Next, it is sufficient to prove that on the event $E$, $\sigma_j^2>0$, which implies $I_E=0$, and so $\pr(\varpi_j=p)=\pr(\bm V\bm u_j^t=p)=0$.
  We denote $\alpha_i=\bm u_idiag(\bm v_i)\bm u_i^t$. Then $\alpha_i>0$, $\bm U^tdiag(\bm V)\bm U=\sum_{i=1}^{\nu_1} \varpi_i \alpha_i \bm v_i^t\bm v_i$. Here, we use the fact that $diag(\bm v_i)\bm u_j^t=0$ for $i\ne j$ due to $\bm v_i\bm u_j^t=0$. So,
 $$\sigma_j^2\ge \bm u_j\bm\Sigma_1\bm u_j^t=\varpi_j(1-\varpi_j)^2 \alpha_j+\varpi_j^2\sum_{i\ne j} \varpi_i\alpha_i>0\;\; \text{when } 0<\varpi_j<1. $$
 Because we have shown that $\varpi_j>0$ a.s., the proof is now completed. $\Box$

\bigskip
{\bf Proof of Theorem \ref{th2.4}.}  Note that (\ref{eqth3.4.1}) and (\ref{eqth3.4.2}) are satisfied. Theorem \ref{th2.4} follows from Theorem \ref{th3.4} immediately. $\Box$

 \end{document}